\documentclass[11pt,a4paper,reqno]{amsart}
\usepackage{preemble-art}

\begin{document}

\title[Glauber-Exclusion dynamics]
      {Glauber-Exclusion dynamics: rapid mixing regime}

      \author[R. Tanaka]{Ryokichi Tanaka}
      %\address{Mathematical Institute, Tohoku University, Sendai 980-8578, JAPAN}
      %\email{ryokichi.tanaka@tohoku.ac.jp}
      \address{Department of Mathematics, 
      Kyoto University, Kyoto 606-8502 JAPAN}
      \email{rtanaka@math.kyoto-u.ac.jp}

	\author[K. Tsunoda]{Kenkichi Tsunoda}
	\address{Faculty of Mathematics,
	Kyushu University, Fukuoka 819-0395, JAPAN}
	\email{tsunoda@math.kyushu-u.ac.jp}

\subjclass[2010]{82C22, 60J27, 82C20}
\keywords{Glauber-Exclusion process, mixing times for Markov chains, interacting particle systems, hydrodynamic limit}
\date{\today}
%82C22 Interacting particle systems in time-dependent statistical mechanics [See also 60K35]
%60J27 Continuous-time Markov processes on discrete state spaces
%60K35 Interacting random processes; statistical mechanics type models; percolation theory [See also 82B43, 82C43]
%82C20 Dynamic lattice systems (kinetic Ising, etc.) and systems on graphs in time-dependent statistical mechanics

\maketitle

\begin{abstract}
We show that for any attractive Glauber-Exclusion process on the one-dimensional lattice of size $N$ with periodic boundary condition, if the corresponding hydrodynamic limit equation has a reaction term with a strictly convex potential, 
then the total-variation mixing time is of order $O(\log N)$.
In particular, the result covers the full high-temperature regime in the original model introduced by De Masi, Ferrari and Lebowitz (1985).
\end{abstract}

\section{Introduction}

The Glauber-Exclusion dynamics is a stochastic spin system which is a superposition of a Glauber dynamics and a speeded-up symmetric simple exclusion process (SSEP) on lattices.
The model was introduced by De Masi, Ferrari and Lebowitz in \cite{DMFL-PRL} and in \cite{DMFL},
where they showed that the hydrodynamic scaling limit yields a reaction-diffusion equation (see also \cite{KOV} and \cite{DMPbook}).
The family of limiting equations possibly interpolates dynamics of different natures in the macroscopic scale where the reaction term admits a single-well or a double-well potential.
A salient feature of this model is that both types can arise even when the dynamics is constrained in the one-dimensional lattice. 
It seems plausible to expect, but still challenging to verify that the time needed for the stochastic model to reach an equilibrium would reflect the corresponding type in the hydrodynamic limit. 
The purpose of this article is to show a rapid mixing property of the underlying stochastic model when it is hinted by the hydrodynamic equation.
More precisely, we consider a general class of Glauber-Exclusion models on the one-dimensional integer lattice of size $N$ with periodic boundary condition,
and show that if the limiting reaction-diffusion equation has a strictly convex single well potential,
then the total-variation mixing time is of order $O(\log N)$.
In particular, this covers the full high-temperature regime in the original model of De Masi {\it et al}.\ \cite{DMFL-PRL, DMFL}.

For any integer $N \ge 1$, let $\Z_N:=\Z/N\Z$ and $\O_N:=\{-1, 1\}^{\Z_N}$.
The generator of Glauber (spin flip) dynamics is defined by
\[
L_G f(\y):=\sum_{x \in \Z_N}c(x, \y)\(f(\y^x)-f(\y)\) \quad \text{for $\y \in \O_N$},
\]
where $\y^x$ is the configuration whose spin is $-\y(x)$ at $x$ and coincides with the spin of $\y$ at the other sites, and $c(x, \y)$ is a positive value called {\it jump rate}.
In the associated continuous-time Markov chain, the transition $\y \to \y^x$ occurs at rate $c(x, \y)$.
We assume that $c(x, \y)$ depends only on $(\y(y))_{|y-x| \le K}$ for some $K \ge 0$ independent of $N$, is periodic, i.e., $c(x, \y)=f(\y_{\cdot + x})$ for some function $f$
where $\y_{\cdot +x}:=(\y(y+x))_{y \in \Z_N}$,
and that there exists a constant $c_0>0$ independent of $N$ such that $c(x, \y) \ge c_0$ for all $(x, \y) \in \Z_N \times \O_N$.
The generator of SSEP is defined by
\[
L_E f(\y):=\frac{1}{2}\sum_{x \in \Z_N}\(f(\y^{x, x+1})-f(\y)\) \quad \text{for $\y \in \O_N$},
\]
where $\y^{x, x+1}$ is the configuration whose spins are obtained by exchanging the spins of $\y$ at $x$ and $x+1$.
Let us consider the sum of generators where the SSEP is speeded-up at rate $N^2$,
\[
L_N:=L_G+N^2 L_E,
\]
and call the associated continuous-time Markov chain $\{\y_t\}_{t \ge 0}$ on $\O_N$ a {\it Glauber-Exclusion} process
\footnote{This is also called a Glauber+Kawasaki process in literature.}.
Further if the jump rate $c(x, \y)$ for the Glauber dynamics is attractive (see the definition in Section \ref{Sec:monotone}),
then we call the corresponding processes an {\it attractive} Glauber-Exclusion process.
Note that the Glauber-Exclusion process is irreducible on $\O_N$ since the jump rate in the Glauber part is positive,
and thus there exists a unique stationary distribution $\pi_N$ on $\O_N$.
We point out that
except for a very special case the process is not reversible as it was observed by Gabrielli {\it et al}.\ \cite[Section 3]{GJLLV} (see Remark \ref{Rem:nonreversible}), and no explicit description of the stationary distribution is available.

For any $\y \in \O_N$, let $\Pr_\y$ be the distribution of the process $\{\y_t\}_{t \ge 0}$ starting from $\y$.
We define the {\it total-variation mixing time} for each $0<\d<1$,
\[
t^N_\mix(\d):=\inf\Big\{t \ge 0 \ : \ \max_{\y \in \O_N}\|\Pr_\y(\y_t \in \cdot \ )-\pi_N\|_{\TV} \le \d\Big\},
\]
where the total-variation distance is defined by
\[
\|\m-\n\|_\TV:=\max_{A \subset \O_N}|\m(A)-\n(A)|=\frac{1}{2}\sum_{\y \in \O_N}|\m(\y)-\n(\y)|,
\]
for any pair of probability distributions $\m, \n$ on $\O_N$.
For $\rho \in [-1, 1]$, 
let us define $R(\rho):=\E_{\n_\rho}[-2\y(0)c(0, \y)]$
where $\n_\rho$ is the product of Bernoulli measures on $\{-1, 1\}$ with mean $\rho$.
The function $R: [-1, 1] \to \R$ is referred to as the reaction function and $R(\rho)$
 as it appears in \eqref{Eq:HDL} is referred to as the reaction term.
A function $V$ is called a {\it potential} if $V'(\rho)=-R(\rho)$ for $\rho \in [-1, 1]$.
We say that $V$ is strictly convex if $V''(\rho)>0$ on $[-1, 1]$.

\begin{theorem}\label{Thm:mixing}
For any attractive Glauber-Exclusion process on $\Z_N$,
if the corresponding reaction function
$R$ admits a strictly convex potential function $V$, 
then there exists a constant $C$ such that for all $0<\d<1$, we have
\[
t^N_\mix(\d) \le \frac{1}{\k}\log \frac{N}{\d}+C,
\]
for all large enough $N$, where $\k:=\min_{\rho \in [-1, 1]}V''(\rho)>0$.
\end{theorem}

For example, let us consider the jump rate in the Glauber part introduced by De Masi {\it et al}.\ \cite{DMFL-PRL, DMFL}:
for $0 \le \g<1$,
\begin{equation}\label{Eq:DMFL-c}
c(x, \y):=1-\g \y(x)(\y(x+1)+\y(x-1))+\g^2 \y(x+1)\y(x-1),
\end{equation}
for $(x, \y) \in \Z_N \times \O_N$.
This jump rate is positive and attractive for any $0 \le \g <1$.
The reaction term
$R(\rho)=\E_{\n_\rho}[-2\y(x)c(x, \y)]=-2(1-2\g)\rho-2\g^2\rho^3$
has a strictly convex potential function if $0 \le \g<1/2$. 
Therefore this range $0 \le \g <1/2$ may be considered as the high temperature regime analogous to the mean-field setting.
Theorem \ref{Thm:mixing} implies that
for every $0\le \g < 1/2$ there exists a constant $C$ such that for all $0<\d<1$, we have
\[
t^N_\mix(\d) \le \frac{1}{2(1-2\g)}\log \frac{N}{\d}+C,
\]
for all large enough $N$.
However, we believe that the constant in front of $\log N$ in Theorem \ref{Thm:mixing}
would not be optimal and leave it open to find an exact constant which might lead to a cutoff particularly for the special case 
\eqref{Eq:DMFL-c} by De Masi {\it et al}.
For background and discussions on the cutoff phenomenon, see \cite[Chapter 18]{LevinPeresBook}.

Let us recall the result from the hydrodynamic limit which has been obtained by Kipnis, Olla and Varadhan \cite[Appendix]{KOV} in the case of finite volume $\Z_N$ and by De Masi {\it et al}.\ \cite{DMFL} and \cite[Chapter VI]{DMPbook} in the case of infinite volume $\Z$.
For background information in the theory of hydrodynamic limit, see \cite{KipnisLandim}; for recent progresses, see \cite{BBP} and references therein.
For each fixed time $t\ge 0$, as $N \to \infty$ under an appropriate condition of convergence of the initial data, the density field 
\[
\m^N_t(dx)=\frac{1}{N}\sum_{x \in \Z_N}\y_t(x)\d_{\frac{x}{N}}(dx) \quad \text{on $\R/\Z$},
\]
converges weakly to a unique weak solution of the reaction-diffusion equation:
	\begin{equation}\label{Eq:HDL}
	\frac{\partial \rho}{\partial t}=\frac{1}{2}\frac{\partial^2\rho }{\partial x^2} +R(\rho),
	\end{equation}
where $R(\rho)$ is the reaction term defined just before the statement of Theorem \ref{Thm:mixing}.
The following is a reason why the reaction term is obtained by the average under the Bernoulli measure:
for a fixed time $T>0$
since the SSEP is speeded-up by $N^2$ whereas the spin-flip dynamics occurs slowly $O(N)$ in a unit time,
in each (macroscopically small) box of size $\e N$ for $\e>0$
the distribution at time $T$ is approximately the Bernoulli measure whose parameter is the empirical average of spins in the box.
There one should be able to replace the distribution by the Bernoulli at each fixed time $T$ --- the proof of deriving \eqref{Eq:HDL} is basically devoted to verify this heuristic rigorously.
Let us point out that this has been verified only for a fixed time $T$, which is not enough to understand the full dynamics up to the equilibrium.
One observes that a combination of (microscopically) slow reaction with a fast stirring creates a non-local interaction among spins
even though the Glauber updates depend only on configurations in bounded ranges (e.g., nearest-neighbors).
Because of this particular nature of the dynamics our analysis would not simply boil down to either the Glauber dynamics or the SSEP,
where much sharper results have been established for the total-variation mixing time in \cite{LubetzkySly} and \cite{Lacoin2016, Lacoin2017}, respectively.
Concerning the Glauber-Exclusion process, very little is known for the stationary distribution, which has attracted intensive studies 
recently (e.g., \cite{FLT}).
Although it is not always necessary to have full understanding of the stationary distribution (e.g., a sharp result on the mixing time has been established for a large class of noisy voter models in \cite{CoxPeresSteif}),
it had been unknown that the mixing time was actually $O(\log N)$ in the present setting.
Given the result on the hydrodynamic limit, our result would not be valid beyond the strict convexity condition on potentials:
even for the particular example by De Masi {\it et al}.\ \eqref{Eq:DMFL-c},
it is tempting to determine the exact order of total-variation mixing time in the range $1/2 \le \g<1$.

Our proof is based on an enhancement of the above replacement which is now valid in time range up to $N^\e$ for a small enough $\e>0$ with a uniform error control.
We call it a {\it replacement lemma} (see the precise statement in Lemma \ref{Lem:replacement});
the proof uses a classical coupling between a SSEP and independent simple random walks
and it sharpens estimates in \cite{DMFL}. 
The idea has a lot of similarities to the ``$v$-function'' method in \cite[Chapter IX]{DMPbook},
however our proof is self-contained and enables us to expose an explicit constant factor in the upper bound if not optimal.
We use this estimate under the standard monotone coupling process
in order to control error terms.
After obtaining a recursive estimate with a fine error control,
we find a coalescing time of the coupling process and this leads an upper bound of the total-variation mixing time (Section \ref{Sec:main}).
In the entire part of this paper, we focus on the one-dimensional lattice with periodic boundary condition, 
but this is mainly for simplicity of notations; the methods we employ are valid in the higher dimensional lattices.

Let us give some other explicit examples to which our result applies.

\begin{example}[The De Masi-Ferrari-Lebowitz model and its variants]
	Brassesco {\it et al}.\ have considered the following jump rate with external field \cite{BPSVm} (compare with \cite{BPSVexp}):
	\[
	c_\m(x,\y):=c(x, \y)-\frac{\m}{2}\y(x),
	\]
	where $c(x, \y)$ is the jump rate of De Masi {\it et al}.\ \eqref{Eq:DMFL-c}, and $\m$ is a real value. 
	If
	$0\le \m<2(1-\g)^2$,
	then $c_\m(x, \y)>0$.
	Note that introducing the parameter $\m$ breaks the flip symmetry $\y \leftrightarrow -\y$ of the stationary distribution and still keeps the jump rate attractive.
	The reaction term has the following form:
	$R(\rho)=\E_{\n_\rho}[-2\y(x)c(x, \y)]=-2(1-2\g)\rho-2\g^2\rho^3+\m$,
	which has a strictly convex potential function if $0 \le \g<1/2$.
	(It should be noted however that the main purpose to introduce $\m$ is to understand the dynamics when $\g>1/2$ in \cite{BPSVm}.)
	Concerning a more general class of processes in the one-dimensional lattice, see \cite[Section 4]{Ferrari}.
	\footnote{In the same paper, there is a typographical error where the Glauber-Exclusion process with the jump rate \eqref{Eq:DMFL-c} is analyzed; it is in fact $\g <1/3$ instead of $\g \le 1/2$ \cite[p.1528, line 22]{Ferrari}.}
\end{example}

\begin{remark}
	For any $d \ge 1$, one is able to generalize the process on the $d$-dimensional lattice by setting the jump rate
	\[
	c(x, \y)=\frac{1}{d}\sum_{i=1}^d\(1-\g \y(x)\(\y(x-e_i)+\y(x+e_i)\)+\g^2 \y(x-e_i)\y(x+e_i)\),
	\]
	where $\{e_i\}_{i=1, \dots, d}$ is the standard basis in $\R^d$.
\end{remark}

\begin{example}[The Chafee-Infante model]
This model has been analyzed in \cite[Section 8]{FLT}.
Let us fix $a_0, a_1, a_2>0$, and define
\begin{align*}
c(x, \y) :=a_0\1_{\{\y(x-1)=\y(x+1)\neq \y(x)\}}+a_1 \1_{\{\y(x-1)=\y(x)=\y(x+1)\}}+a_2 \1_{\{\y(x-1)\neq \y(x+1)\}}.
\end{align*}
Note that this jump rate is positive and attractive if $a_0 \ge a_2 \ge a_1$.
A computation yields
\begin{align*}
\E_{\n_\rho}[-2\y(x)c(x, \y)]=\frac{1}{2}\((a_0-3a_1-2a_2)\rho-(a_0+a_1-2a_2)\rho^3\).
\end{align*}
Fix $a, b>0$.
If we define 
$a_1:=a$, $a_2:= a+2b$ and
$a_0:=a+8b$, then the jump rate is positive and monotone and we have
\[
a_0-3 a_1-2a_2=4(b-a)
\qquad
\text{and}
\qquad
a_0+a_1-2a_2=4b.
\]
Hence
\[
R(\rho)=\E_{\n_\rho}[-2\y(x)c(x, \y)]=2(b-a)\rho-2 b \rho^3,
\]
where $R(\rho)$ has a strictly convex potential function if $b < a$.
\end{example}

\begin{remark}\label{Rem:nonreversible}
In general, an explicit form of the stationary distribution $\pi_N$ is not available even in the setting which we are focusing on in this article. 
Gabrielli {\it et al}.\ have shown that the Glauber-Exclusion process on periodic integer lattices is reversible if and only if the jump rate has the form
\[
c(x, \y)=(a_1+a_2 \y(x))h(x, \y),
\]
where $h(x, \y)$ is independent of $\y(x)$ with some constants $a_1, a_2$, and
in this particular case $\pi_N$ is a Bernoulli measure \cite[Section 3]{GJLLV}.
\end{remark}

The rest of this article is organized as follows:
in Section \ref{Sec:main} we prove Theorem \ref{Thm:mixing} provided Lemma \ref{Lem:nonlinear} whose full proof is deferred until subsequent sections,
in Section \ref{Sec:proof of lemma nonlinear}
we prove Lemma \ref{Lem:nonlinear} by using a key replacement lemma Lemma \ref{Lem:replacement},
in Section \ref{Sec:taming} we prove Lemma \ref{Lem:replacement} and complete the proof of upper bound for the mixing time,
and in Appendix \ref{Sec:A} and \ref{Sec:misc} we give the proofs of technical lemmas which we use in the proof of Lemma \ref{Lem:replacement} in Section \ref{Sec:taming}.

{\bf Notation}:
Throughout the article, we use $C, C', C'', \dots$ to denote numerical constants whose exact values may change from line to line, and $C_\e, C'_\e, C''_\e, \dots$ to indicate its dependency only on $\e$ for a parameter $\e$.
For a set $A$, we denote the indicator function of $A$ and the cardinality of $A$ by $\1_A$ and $|A|$, respectively.
When $A$ is a Lebesgue measurable subset of $\R$, we also denote the normalized Lebesgue measure of $A$ in $\R$ by $|A|$, for which we believe there is no danger of confusion from the context.
For a real-valued function $f(N)$ in positive integers $N$, we write $f(N)=O(N)$ if there exists a constant $C\ge0$ such that 
$|f(N)| \le C N$ for all large enough $N$, 
and $f(N)=o(1)$ if $|f(N)| \to 0$ as $N \to \infty$.
For any real number $x$, we denote by $\lfloor x\rfloor$ the largest integer at most $x$ and by $\lceil x\rceil$ the smallest integer at least $x$.

\section{Proof of Theorem \ref{Thm:mixing}}\label{Sec:main}

\subsection{Monotone coupling}\label{Sec:monotone}

We define a partial order in $\O_N$ by setting
\[
\y \ge \z \iff \y(x) \ge \z(x) \quad \text{for any $x \in \Z_N$}.
\]
A jump rate $c(x, \y)$ $(x \in \Z_N, \y \in \O_N)$ is called
{\it attractive} if $\y \ge \z$ implies that
\[
c(x, \y) \le c(x, \z) \text{\ for $\y(x)=\z(x)=1$},
\quad
\text{and}
\quad
c(x, \y) \ge c(x, \z) \text{\ for $\y(x)=\z(x)=-1$}.
\]
An attractive jump rate yields a monotone coupling of the associated chains $(\y_t^+, \y_t^-)$, i.e., a coupling where the order of configurations are kept preserved $\y_t^+ \ge \y_t^-$ for all $t \ge 0$ almost surely if they have started with $\y_0^+\ge \y_0^- $ \cite[Theorem 4.11, p.143]{LiggettBook}.
We construct a monotone coupling $(\y_t^+, \y_t^-)$ for the Glauber-Exclusion process based on this coupling, for which we provide the explicit form for the sake of convenience:
given $\y \ge \z$ (otherwise we define the jump rate $0$ from $(\y, \z)$),
for each $x \in \Z_N$,
\[
(\y, \z) \mapsto (\y^{x, x+1}, \z^{x, x+1}) \quad \text{at rate $\frac{N^2}{2}$},
\]
and the states of $(\y, \z)$ at $x$ are updated according to the monotone coupling of Glauber dynamics with the same (attractive) jump rate, namely, we change the states of $(\y, \z)$ at $x$ by the following rates;
\begin{align*}
(\y(x), \z(x))=(-1, -1) &\implies 
\begin{cases}
(1, 1)  & \ \ \text{with rate $c(x, \z)$}\\
(1, -1) & \ \ \text{with rate $c(x, \y)-c(x, \z)$},
\end{cases}\\
(\y(x), \z(x))=(1, -1) &\implies 
\begin{cases}
(-1, -1) & \text{with rate $c(x, \y)$}\\
(1, 1) & \text{with rate $c(x, \z)$},
\end{cases}\\
(\y(x), \z(x))=(1, 1) &\implies 
\begin{cases}
(-1, -1) & \text{with rate $c(x, \y)$}\\
(1, -1) & \text{with rate $c(x, \z)-c(x, \y)$}.
\end{cases}
\end{align*}
The joint process 
$\{(\y_t^+, \y_t^-)\}_{t\ge 0}$ yields a Markovian coupling, whose distribution we denote by the same symbol $\Pr$ as long as there is no danger of confusion.
If we define
\[
\t:=\inf\{t \ge 0 \ : \ \y^+_t =\y^-_t\},
\]
then it holds that $\t < \infty$ and $\y_t^+=\y_t^-$ for all $t \ge \t$ almost surely if $\y^+_0 \ge \y^-_0$ since the chain is irreducible.
Let $\1$ and $-\1$ be the configurations whose states are all $1$ and all $-1$, respectively.

\begin{lemma}\label{Lem:monotone-coupling}
For any attractive Glauber-Exclusion process,
if the corresponding reaction term $R(\rho)$ has a potential $V(\rho)$, i.e., $R(\rho)=-V'(\rho)$ such that 
\[
\k:=\min_{\rho \in [-1, 1]}V''(\rho)>0,
\]
then for the monotone coupling $(\y^+_t, \y^-_t)$ with the initial configurations $(\1, -\1)$,
there exists a constant $C>0$ such that for any $0<\d<1$ and
\[
t_{N, \d}:=\frac{1}{\k}\log \frac{N}{\d}+C, 
\]
we have
$\Pr(\t>t_{N, \d}) \le \d$ for all large enough $N$.
\end{lemma}

We show Lemma \ref{Lem:monotone-coupling} in Section \ref{Sec:proof of Lem coupling}.

\proof[Proof of Theorem \ref{Thm:mixing}]
Let $\y, \z \in \O_N$  be two arbitrary configurations (not necessarily $\y \ge \z$).
We consider the monotone coupling $\{(\y_t^+, \y_t^-)\}_{t\ge 0}$ with initial configurations $(\1, -\1)$.
Note that $\1 \ge \y, \z \ge -\1$;
enlarging the probability space,
we construct a coupling $\Pr$ (denoted by the same symbol) among four copies of original chains $\y^+_t, \y^-_t, \y_t, \z_t$ such that 
\[
\y^+_t \ge \y_t, \z_t \ge \y^-_t \quad \text{for all $t \ge 0$ almost surely},
\]
with $\y_0=\y$ and $\z_0=\z$, and further
any two of $\y_t^+$, $\y_t$ and $\y_t^-$ (resp.\  $\y_t^+$, $\z_t$ and $\y_t^-$) are monotone coupling.
Then $\y^+_t=\y^-_t$ implies that $\y_t=\z_t$ almost surely, and thus
letting 
\[
\t:=\inf\{t \ge 0 \ : \ \y_t^+=\y_t^-\},
\]
we have that
\[
\|\Pr_\y(\y_t \in \cdot \ )-\Pr_\z(\z_t \in \cdot \ )\|_\TV \le \Pr(\y_t \neq \z_t) \le \Pr(\y_t^+ \neq \y_t^-) \le \Pr(\t>t).
\]
By Lemma \ref{Lem:monotone-coupling} for any $0<\d<1$
if
$t_{N, \d}=(1/\k)\log (N/\d)+C$ for some constant $C$,
then
$\Pr(\t>t_{N, \d}) \le \d$.
Noting that
\[
\max_{\y \in \O_N}\|\Pr_\y(\y_t \in \cdot \ ) - \pi_N\|_\TV \le \max_{\y, \z \in \O_N}\|\Pr_\y(\y_t \in \cdot \ )-\Pr_\z(\z_t \in \cdot \ )\|_\TV,
\]
we obtain 
\[
t^N_\mix(\d) \le \frac{1}{\k}\log \frac{N}{\d}+C,
\]
for all large enough $N$ and conclude the proof.
\qed

\subsection{The upper bound}\label{Sec:proof of Lem coupling}

For a continuous-time Markov chain $\{\y_t\}_{t\ge 0}$ generated by $L_N$ on the state space $\O_N=\{-1, 1\}^{\Z_N}$,
we define the {\it normalized magnetization} by
\[
S(\y_t):=\frac{1}{N}\sum_{x \in \Z_N}\y_t(x).
\]

We say that $f:\{-1, 1\}^{\Z_N} \to \R$ is a {\it local function} if there exists an integer $0 \le K <N$ independent of $N$
and a subset $I$ in $\Z_N$ with $|I|\le K$
such that $f(\y)$ depends only on $\{\y_x\}_{x \in I}$ for any $\y \in \{-1, 1\}^{\Z_N}$.
Let us denote by $\supp f$ the smallest such subset $I$ for a local function $f$.
For example, $f(\y):=\y(-1)\y(0)\y(1)$ for $\y \in \{-1, 1\}^{\Z_N}$ is a local function with $\supp f=\{-1, 0, 1\}$.
Let $f_\emptyset(\y)=1$.
For a non-empty subset $I$ in $\Z_N$,
let $f_I(\y):=\prod_{x \in I}\y(x)$
and call it an {\it elementary local function}.
Note that any local function $f$ is a linear combination of elementary local functions:
\[
f(\y)=\sum_{I \subset \supp f}a_I f_I(\y) \quad \text{for $\y \in \O_N$},
\]
where the summation is over all subsets in $\supp f$ and $a_I \in \R$
and this follows from the Fourier expansion of $f$ on $\{-1, 1\}^K$.

For $\y \in \O_N$ and for $x \in \Z_N$,
denoting
$L_N f_x(\y)$ where $f_x(\y):=\y(x)$ simply by $(L_N\y)(x)$, 
we have
\begin{align*}
	(L_N \y)(x)	%=(L_G\y)(x)+N^2 (L_E \y)(x)
	=-2 \y(x)c(x, \y)+\frac{N^2}{2}\(\y(x+1)+\y(x-1)-2\y(x)\),
\end{align*}
and thus
\[
L_N S(\y)=\frac{1}{N}\sum_{x \in \Z_N}-2 \y(x)c(x, \y).
\]
Since for a given jump rate $c(x, \y)$, the function $\y(0)c(0, \y)$ is local,
expanding it by elementary local functions in the way above, we have
\[
2\y(0)c(0, \y)=\sum_{y\in\mathbb Z_N}a_y f_y(\y)+g(\y) \quad \text{where $g(\y):=\sum_{\text{$|I|\neq1$}}a_I f_I$}.
\]
Note that $\eta(x)c(x,\eta)=\eta_{\cdot+x}(0)c(0,\eta_{\cdot+x})$ and that $\frac1N\sum_{y\in\mathbb Z_N} f_y(\eta_{\cdot+x})=S(\eta)$
to obtain the following decomposition:
\[
L_N S(\y)=-\a S(\y)-F(\y) \quad \text{where $\a=\sum_{y \in \Z_N} a_y$ and $F(\y):=\frac{1}{N}\sum_{x \in \Z_N}g(\y_{\cdot +x})$}.
\]
Note that the reaction term $R(\rho)$ is obtained from $L_N S(\y)$ by
\[
R(\rho)=\E_{\n_\rho}[L_N S(\y)].
\]
Letting 
$G(\rho):=\E_{\n_\rho}F(\y)$,
one has 
$R(\rho)=-\a \rho-G(\rho)$.
If $V(\rho)$ is a potential function of $R(\rho)$, i.e., $V'(\rho)=-R(\rho)$,
then we have that $V''(\rho)=\a+G'(\rho)$.

For a pair of configurations $\y^+_0, \y^-_0$ such that $\y^+_0 \ge \y^-_0$, we run the monotone coupling $\{(\y_t^+, \y_t^-)\}_{t \ge 0}$ and denote by $\Fc_t:=\s\(\{\y_s^+, \y_s^-\}_{0 \le s \le t}\)$ 
the $\s$-algebra generated by the monotone coupling process up to time $t$ for $t \ge 0$.
Letting 
\[
\x_t:=S(\y_t^+)-S(\y_t^-)
\]
for $t \ge 0$,
we have that
\begin{align}\label{Eq:xi}
	\E [\x_{t+T}\mid \Fc_t]=\x_t-\a\int_t^{t+T} \E[\x_s\mid \Fc_t]\,ds- \int_t^{t+T} \E[F(\y_s^+)-F(\y_s^-)\mid \Fc_t]\,ds,
\end{align}
where $\E[\,\cdot \mid \Fc_t]$ stands for the conditional expectation with respect to the $\s$-algebra $\Fc_t$.
For each $x \in \Z_N$, 
we define the {\it local average} of $\y \in \O_N$ around $x$ up to time $t$ by
letting $\{w(t)\}_{t \ge 0}$ be the simple random walk (SRW) starting from $x$ with rate $1$,
\[
\F_x(\y, t):=\Eb_x \y(w(t)) \quad \text{for $\y \in \O_N$ and $t \ge 0$},
\]
where $\Eb_x$ denotes the expectation with respect to the distribution of $\{w(t)\}_{t \ge 0}$ such that $w(0)=x$.

\begin{lemma}[Main replacement lemma]\label{Lem:nonlinear}
For all small enough $\e>0$,
there exist constants $C>0$ and $T_0$ such that 
for all $T_0/N^2 \le t \le 1/N^{1-\e}$ and all $\wt \y=(\y^+, \y^-)$ with $\y^+ \ge \y^-$, we have
\begin{align*}
	&\left|\E_{\wt \y}[F(\y_{t}^+)-F(\y_{t}^-)]-\frac{1}{N}\sum_{x \in \Z_N}[G(\F_x(\y^+, N^2 t))-G(\F_x(\y^-, N^2 t))]\right|\\
	&\qquad \qquad \qquad \qquad \qquad \qquad \qquad \qquad \le \frac{C}{t^{\e/4} N^{2\e}}\x_0
	+C \int_0^t \E \x_s\,ds+C \exp(-(N^{2}t)^{\e/8}).
\end{align*}
\end{lemma}

We show Lemma \ref{Lem:nonlinear} in Section \ref{Sec:proof of lemma nonlinear}
and now prove Lemma \ref{Lem:monotone-coupling}.
We note that the proof of Lemma \ref{Lem:replacement} reveals that $0<\e<1/28$ is sufficient for the argument.
The proof of Lemma \ref{Lem:monotone-coupling} is based on showing a strong contraction (see \eqref{Eq:xi3}) between two coupled magnetization chains.
This is achieved by dealing with the function $F$, which is a main `non-linear' term in this discrete setting.

\proof[Proof of Lemma \ref{Lem:monotone-coupling}]
Fix $0<\e<1/28$ and let $T_\ast:=1/N^{1-\e}$. We first claim that for all $t\ge0$,
\begin{equation}\label{Eq:xi3}
\E \x_{t+T_\ast} \le \(1-\frac{\k}{N^{1-\e}}+\frac{C_{\e}}{N^{1+\e^2}}\)\E \x_t+\frac{C'_\e}{N^2}.
\end{equation}
We prove this inequality later and let us conclude Lemma \ref{Lem:monotone-coupling} by invoking this estimate.
For any $0<\d<1$, if we define 
\begin{equation*}\label{Eq:t}
t_{N, \d}:=\frac{1}{\k}\(\log N+\log \frac{4}{\d}+C_\e'\),
\end{equation*}
then applying \eqref{Eq:xi3} inductively $\lfloor t_{N, \d}/T_\ast \rfloor$-times and using the monotone coupling $\E \x_t \le \E \x_s$ for $ t\ge s$ yield
\begin{align*}\label{Eq:xi4}
\E \x_{t_{N, \d}}
&\le  \(1-\frac{\k}{N^{1-\e}}+\frac{C_{\e}}{N^{1+\e^2}}\)^{\lfloor t_{N, \d}/T_\ast \rfloor}\E \x_0+\frac{C'_\e}{N^2}\sum_{n=0}^{\lfloor t_{N, \d}/T_\ast \rfloor-1} \(1-\frac{\k}{N^{1-\e}}+\frac{C_{\e}}{N^{1+\e^2}}\)^n\\
&\le \x_0 \exp\(-\k t_{N, \d}+C'_\e\)+N^{1-\e}\frac{C''_\e}{N^2}
%\le \x_0 \exp\(-\(\log N+H\)+C'_\e\)+\frac{C''_\e}{N^{1+\e}}
\le \frac{\d}{2N}+\frac{C''_\e}{N^{1+\e}},
\end{align*}
for all large enough $N$, where we have used $\x_0 \le 2$ in the last inequality.
Recalling that 
$\t:=\inf\{t \ge 0 \ : \x_t =0\}$,
we obtain by the Markov inequality,
\begin{equation*}\label{Eq:xi5}
\Pr\big(\t >t_{N, \d} \big) = \Pr\big(\x_{t_{N, \d}} \ge 1/N\big)\le  N \cdot \E \x_{t_{N, \d}}  \le \frac{\d}{2}+\frac{C''_\e}{N^\e}<\d,
\end{equation*}
for all large enough $N$, as required.

Let us turn to the claim \eqref{Eq:xi3}.
We run the monotone coupling with the initial state $(\1, -\1)$.
Note that for all $s, t>0$,
\begin{equation} \label{Eq:S}
\left|\int_s^{s+t}\E[\x_u \mid \Fc_s]\,du-t\x_s\right| \le C t^2.
\end{equation}
Indeed, since the SSEP dynamics given by $L_E$ preserves $S(\y)$, 
the difference 
$|\x_u-\x_s|$ is at most $2/N$ times the number of Glauber updates within the time interval $[s, u]$ for $0 \le s<u$,
which is stochastically dominated by the Poisson random variable with intensity $O(N (u-s))$.
Integrating in $u$,
\[
\E[\x_u-\x_s \mid \Fc_s] \le \E[|\x_u-\x_s|\mid \Fc_s] \le C(u-s),
\]
we obtain \eqref{Eq:S}.
For each $t \ge 0$, we integrate
\[
\E[F(\y_{s}^+)-F(\y_{s}^-) \mid \Fc_t] \quad \text{for $t \le s \le t+T_\ast$},
\]
where we decompose the time interval for the integration into two parts; for the first part $t \le s \le t+T_0/N^2$, we use the following elementary estimate
\[
\left| \int_t^{t+\frac{T_0}{N^2}} \E[F(\y_s^+)-F(\y_s^-) \mid \Fc_t]\,ds\right|\le \frac{2\|g\|_\infty T_0}{N^2},
\]
(for the definition of $g$, see the third paragraph in Section \ref{Sec:proof of Lem coupling}) and for the second part $t+T_0/N^2 \le s \le T_\ast$,
we will apply to the Markov property and Lemma \ref{Lem:nonlinear}.
By \eqref{Eq:S}, it follows that
\[
\int_t^{t+s} \E[\xi_u\mid \Fc_t]\,du \le s \xi_t +C s^2,
\]
which we integrate over $t \le s \le t+T_\ast$ (and we note that this factor comes from the second term in the right hand side in the inequality in Lemma \ref{Lem:nonlinear}).
Direct computations yield
\[
\int_0^{T_\ast} \(s+\frac{1}{s^{\e/4} N^{2\e}}\)\,ds =\frac{1}{2}T_\ast^2+ \frac{1}{(1-\e/4)N^{2 \e}}T_\ast^{1-\e/4} \le \frac{C_\e}{N^{1+\e^2}},
\]
(which is a factor in front of $\x_t$; we have used $(1-\e)(1-\e/4)+2\e \ge 1+\e^2$),
and
\[
\int_0^{T_\ast} \(s^2+\exp\(-(N^{2}s)^{\e/8}\)\)\,ds\le \frac{1}{3}T_\ast^3+\frac{1}{N^2}\int_0^\infty \exp\(-s^{\e/8}\)\,ds\le \frac{C'_\e}{N^2},
\]
(which is a remaining error term),
combining above estimates, we have that 
\begin{align}
\int_t^{t+T_\ast} \E[F(\y_{s}^+)-F(\y_{s}^-) \mid \Fc_t]\,ds
&\ge -\frac{C''_\e}{N^2}-\frac{C_\e}{N^{1+\e^2}}\x_t\notag\\
&+\int_0^{T_\ast}\frac{1}{N}\sum_{x \in \Z_N}[G(\F_x(\y_t^+, N^2 s))-G(\F_x(\y_t^-, N^2 s))]\,ds, \label{Eq:4}
\end{align}
where we have used $\x_t \ge 0$, the Markov property and Lemma \ref{Lem:nonlinear}.
For any pairs of configurations with $\y^+ \ge \y^-$,
it holds that for $s \ge 0$,
\[
\F_x(\y^+, N^2 s) \ge \F_x(\y^-, N^2 s).
\]
Recall that a potential function $V(\rho)$ of $R(\rho)$ satisfies that $V''(\rho)=\a+G'(\rho)$.
Letting $\k:=\min_{\rho \in [-1, 1]}V''(\rho)>0$,
we have that by the Taylor theorem
for $\rho^+ \ge \rho^-$, 
\begin{align*}
G(\rho^+)-G(\rho^-)=-\a(\rho^+-\rho^-)+(V'(\rho^+)-V'(\rho^-))
\ge (-\a+\k)(\rho^+-\rho^-).
\end{align*}
Plugging into this inequality $\F_x(\y_t, N^2 s)$ and noting that
\[
\frac{1}{N}\sum_{x \in \Z_N}\F_x(\y_t, N^2 s)=\frac{1}{N}\sum_{x \in \Z_N}\Eb_x \y_t(w(N^2 s))=S(\y_t),
\]
we obtain by averaging over all sites $x \in \Z_N$, for all $s \ge 0$,
\begin{align}\label{Eq:5}
\frac{1}{N}\sum_{x \in \Z_N}[G(\F_x(\y_t^+, N^2 s))-G(\F_x(\y_t^-, N^2 s))] \ge (-\a+\k)(S(\y_t^+)-S(\y_t^-)).
\end{align}
For all large enough $N$, taking expectation in \eqref{Eq:4}, we obtain from \eqref{Eq:xi} and \eqref{Eq:5} for all $t \ge 0$,
\begin{equation}\label{Eq:xi2}
\E[\x_{t+T_\ast} \mid \Fc_t]\le \(1+\frac{C_{\e}}{N^{1+\e^2}}\)\x_t -\a \int_t^{t+T_\ast}\E[\x_s \mid \Fc_t]\,ds+(\a-\k)T_\ast  \x_t+\frac{C'_\e}{N^2}.
\end{equation}
By the monotone coupling, if the initial state is $(\1, -\1)$, then $\E \x_t \le \E \x_s$ for $ t\ge s$,
whence taking expectation in \eqref{Eq:xi2} yields
\[
\(1+\frac{\a}{N^{1-\e}}\)\E\x_{t+T_\ast} \le \(1+\frac{C}{N^{1+\e^2}}\)\E \x_t+(\a-\k)T_\ast \E \x_t+\frac{C_\e'}{N^2}.
\]
Since $T_\ast=1/N^{1-\e}$, multiplying the inverse of $1+\a/N^{1-\e}$ in both sides yields \eqref{Eq:xi3}, which completes the proof of Lemma \ref{Lem:monotone-coupling}.
\qed

\section{Proof of Lemma \ref{Lem:nonlinear}}\label{Sec:proof of lemma nonlinear}

For any positive integer $K \ge 1$ independent of $N$ such that $K \le N$,
let us fix $K$-distinct sites $\{x_1, \dots, x_K\}$ in $\Z_N$.
We consider an SSEP generated by $L_E$ with particles marked by labels $\{1, \dots, K\}$ as a stirring dynamics, i.e., a $K$-marked SSEP $\{z(t)\}_{t \ge 0}$ with rate $1$ starting from 
$z_i(0)=x_i$ 
where $z(t)=(z_i(t))_{i=1, \dots, K} \in (\Z_N)^K$,
and 
for each bond $\{x, x+1\}$ in $\Z_N$, an exchange occurs at random times which are independent Poisson point processes with intensity $1/2$.
Note that for the $N$-marked SSEP $\wt z(t)=(\wt z_x(t))_{x \in \Z_N}$ with rate $1$ starting from $\wt z_x(0)=x$ for $x \in \Z_N$,
we naturally identify $z(t)$ with a subset of $\wt z(t)$ such that $z_{i}(t)=\wt z_{x_i}(t)$ for all $i=1, \dots, K$ and for all $t \ge 0$ almost surely.
Let us define
\[
\y^{z(t)}(x):=
\begin{cases}
\y(z_i(t))	& \text{if $x=x_i$ for $i=1, \dots, K$},\\
\y(x)		& \text{if $x \notin \{x_1, \dots, x_K\}$},
\end{cases}
\]
for $x \in \Z_N$ and $\y \in \O_N$.
This process is constant for all $x\notin \{x_1,\ldots, x_K\}$ and records the values observed by $z_i(t)$ on $\eta$ at time $t$ and each $x_i$ for $i=1,\ldots,K$.
%This process amounts to keep track of the values initially at $x_i$ and they are observed by $z_i(t)$ at time $t$.

\begin{lemma}\label{Lem:integration-by-parts}
	For any local function $f:\O_N \to \R$,
	let
	\[
	F(\y):=\frac{1}{N}\sum_{x \in \Z_N}f(\y_{\,\cdot \,+x}) \quad \text{for $\y \in \O_N$},
	\]
	and $\{z(t)\}_{t \ge 0}$ be the $N$-marked SSEP with rate $1$ starting from each site in $\Z_N$.
	Then there exists a constant $C>0$ depending only on $\|f\|_\infty=\max_{\y \in \O_N}|f(\y)|$, the size of the support $|\supp f|$ and the jump rate such that the following holds:
	for any $t >0$ and for any monotone coupling $(\y^+_t, \y^-_t)$ with the initial configurations $\wt \y:=(\y^+, \y^-)$ where $\y^+ \ge \y^-$,
	\begin{align*}
		\Big|\E_{\wt \y} [F(\y^+_t)-F(\y^-_t)]-\Eb [F(\y^{+z(N^2 t)})-F(\y^{-z(N^2 t)})]\Big|
		\le C\int_0^t \E(S(\y_s^+)-S(\y_s^-))\,ds,
	\end{align*}
	where $\y^{+z(N^2 t)}$ and $\y^{-z(N^2 t)}$ are $\y^{z(N^2 t)}$ for $\y=\y^+$ and $\y^-$, respectively, 
	and $\Eb$ stands for the expectation with respect to $\{z(t)\}_{t \ge 0}$.
\end{lemma}

We defer the proof for a moment and state a technical lemma needed for it.
For any $t>0$,
if we define
\[
\wt f_{t-s}(\y):=\Eb f(\y^{z(N^2(t-s))}) \quad \text{for $s \in [0, t]$ and $\y \in \O_N$},
\]
then
\[
\E \wt f_0(\y_t)-\E \wt f_t(\y)=\int_0^t \E\(L_N+\partial_s\) \wt f_{t-s}(\y_s)\,ds,
\]
where $\partial_s:=\partial/\partial s$. 
Noting that
\[
\partial_s \wt f_{t-s}(\y)=-N^2L_E \wt f_{t-s}(\y),
\]
by the definition of generator $L_N=L_G+N^2 L_E$, we have
\[
(L_N+\partial_s)\wt f_{t-s}(\y_s)=L_N \wt f_{t-s}(\y_s)-N^2 L_E \wt f_{t-s}(\y_s)=L_G \wt f_{t-s}(\y_s),
\]
and $\E \wt f_0(\y_t)=\E f(\y_t)$ 
and $\E \wt f_t(\y)=\Eb f(\y^{z(N^2 t)})$,
we obtain 
\begin{equation}\label{Eq:integration-by-parts}
\E f(\y_t)-\Eb f(\y^{z(N^2 t)})=\int_0^t \E L_G \wt f_{t-s}(\y_s)\,ds \quad \text{for $t > 0$}.
\end{equation}

We use the following lemma to show Lemma \ref{Lem:integration-by-parts}.

\begin{lemma}\label{Lem:local-function}
	For any local function $f:\{-1, 1\}^{\Z_N} \to \R$,
	it holds that
	\[
	f(\y^+)-f(\y^-)=\sum_{x \in \supp f}\(\y^+(x)-\y^-(x)\)\f_{x, f}(\y^+, \y^-),
	\]
	for any $(\y^+, \y^-) \in \O_N^2$,
	where $\f_{x, f}(\y^+, \y^-)$ is a polynomial with respect to $\y^+(y), \y^-(y)$ for $y \in \supp f\setminus \{x\}$ for each $x \in \supp f$.
	Moreover, we have
	\[
	\max_{\y^+, \y^- \in \O_N}|\f_{x, f}(\y^+, \y^-)| \le 2^{|\supp f|}\max_{\y \in \O_N}|f(\y)|.
	\]
\end{lemma}

\proof
For any nonempty subset $I$ in $\Z_N$,
let $\f_I(\y):=\prod_{x \in I}\y(x)$ and $\f_\emptyset(\y)\equiv 1$ for $\y \in \O_N$,
and recall that
any local function $f:\O_N \to \R$ is a linear combination of $\{\f_I\}_{I \subset \supp f}$ with coefficients in $\R$.
Hence it suffices to show the claim for each $\f_I$ with $I \neq \emptyset$.
This follows from the induction on the size of $I$;
indeed,
if $|I|=1$, then $\f_I(\y^+)-\f_I(\y^-)=\y^+(x)-\y^-(x)$ with $I=\{x\}$.
%and $\f_x \equiv 1$.
For $y \notin I$,
\begin{align*}
	\f_{I\cup \{y\}}(\y^+)-\f_{I\cup \{y\}}(\y^-)
	&=\y^+(y) \f_I(\y^+)-\y^-(y) \f_I(\y^-)\\
	&=(\y^+(y)-\y^-(y))\f_I(\y^+)+\y^-(y)\(\f_I(\y^+)-\f_I(\y^-)\).
\end{align*}
This shows that for the base of elementary local functions $\f_I$,
\[
\f_I(\y^+)-\f_I(\y^-)=\sum_{x \in I}(\y^+(x)-\y^-(x))\f_{x, I}(\y^+, \y^-),
\]
where $\f_{x, I}(\y^+, \y^-)$ is a polynomial with respect to $\y^+(y), \y^-(y)$ for $y \in I\setminus\{x\}$ for each $x \in I$, and further $|\f_{x, I}(\y^+, \y^-)| \le 1$.
Since $f=\sum_{I \subset \supp f}a_I \f_I$,
where $a_I:=E_{{\rm supp}  f} [f\varphi_I]$, with $E_{{\rm supp}f}$
the expectation with respect to the product of the $1/2$-Bernoulli distribution on $\{-1,1\}^{{\rm supp}f}$,
$|a_I| \le \|f\|_\infty:=\max_{\y \in \O_N}|f(\y)|$ and we have that
\begin{align*}
f(\y^+)-f(\y^-)
&=\sum_{I \subset \supp f}a_I(\f_I(\y^+)-\f_I(\y^-))\\
&=\sum_{I \subset \supp f}a_I \sum_{x \in \supp f}(\y^+(x)-\y^-(x))\f_{x, I}(\y^+, \y^-)1_{\{x \in I\}}\\
&=\sum_{x \in \supp f}(\y^+(x)-\y^-(x))\sum_{I \subset \supp f}a_I \f_{x, I}(\y^+, \y^-)1_{\{x \in I\}}.
\end{align*}
Letting $\f_{x, f}(\y^+, \y^-):=\sum_{I \subset \supp f}a_I \f_{x, I}(\y^+, \y^-)1_{\{x \in I\}}$,
we obtain $|\f_{x, f}(\y^+, \y^-)| \le 2^{|\supp f|} \|f\|_\infty$, as required.
\qed

\proof[Proof of Lemma \ref{Lem:integration-by-parts}]
For any local function $f$, let $I:=\supp f=\{x_1, \dots, x_{|I|}\}$, and we consider $|I|$-marked SSEP $\{z(t)\}_{t \ge 0}$ with rate $1$ and the initial state $z_i(0)=x_i$ for $i=1, \dots, |I|$.
For any fixed $\y \in \O_N$,
let us write $\wt f_{t-s}(\y):=\Eb f(\y^{z(N^2(t-s))})$, where $\Eb$ stands for the expectation with respect to $z(N^2(t-s))$.
Since for now we fix $t, s$ and $N$,
we just write for the sake of brevity with no danger of confusion %by letting $\Zc:=z(N^2(t-s))$,
\[
\wt f_{t-s}(\y):=\Eb f(\y^\Zc) \quad \text{where $\Zc:=z(N^2(t-s))$},
\]
and $\Zc$ is also considered as a subset in $\Z_N$.
Then we have
\begin{align*}
	L_G \wt f_{t-s}(\y)=\sum_{x \in \Z_N}c(x, \y)\(\wt f_{t-s}(\y^x)-\wt f_{t-s}(\y)\)
	&=\sum_{x \in \Z_N}c(x, \y)\(\Eb f((\y^x)^\Zc)-\Eb f(\y^\Zc)\)\\
	&=\Eb \sum_{x \in \Z_N}c(x, \y)\(f((\y^x)^\Zc)-f(\y^\Zc)\).
\end{align*}
Since $f(\y^\Zc)$ depends only on $\y(y)$ for $y \in \Zc$,
if $x \notin \Zc$, then
$f((\y^x)^\Zc)-f(\y^\Zc)=0$,
and thus
\[
\sum_{x \in \Z_N}c(x, \y)\(f((\y^x)^\Zc)-f(\y^\Zc)\)=\sum_{x \in \Zc}c(x, \y)\(f((\y^x)^\Zc)-f(\y^\Zc)\).
\]

Given $\Zc$ and $x \in \Zc$,
the term $c(x, \y)f(\y^\Zc)$ is a local function depending only on $\y_y$ for $y \in \Zc \cup I_x$,
where $I_x:=\supp c(x, \cdot\,)$, i.e., the support of $c(x, \cdot\,)$.
Hence Lemma \ref{Lem:local-function} implies that
for $\y^+, \y^- \in \O_N$, we have
\[
c(x, \y^+)f(\y^{+ \Zc})-c(x, \y^-)f(\y^{- \Zc})=\sum_{y \in \Zc \cup I_x}(\y^+(y)-\y^-(y))\f_{y, \Zc}(\y^+, \y^-),
\]
where $\f_{y, \Zc}(\y^+, \y^-)$ is a polynomial with respect to $\y^+(z), \y^-(z)$ for $z \in \Zc \cup I_x$,
and moreover, there exists a constant $C$ depending only on the size of the support for $c(x, \y)f(\y^\Zc)$ (which is at most $|\supp f|+|\supp c(x, \cdot )|$),
and $\max_{\y \in \O_N}|c(x, \y)f(\y^\Zc)| \le \|c(0, \cdot\,)\|_\infty \|f \|_\infty$, such that
\[
\max_{\y^+, \y^- \in \O_N}|\f_{y, \Zc}(\y^+, \y^-)|\le C.
\]
This shows that
\[
|c(x, \y^+)f(\y^{+ \Zc})-c(x, \y^-)f(\y^{- \Zc})| \le C\sum_{y \in \Zc \cup I_x}|\y^+(y)-\y^-(y)|,
\]
and we note that the same bound holds for $|c(x, \y^+)f((\y^{+x})^{\Zc})-c(x, \y^-)f((\y^{-x})^{\Zc})|$:
\[
|c(x, \y^+)f((\y^{+x})^\Zc)-c(x, \y^-)f((\y^{-x})^\Zc)| \le C\sum_{y \in \Zc \cup I_x}|\y^+(y)-\y^-(y)|.
\]
Summarizing the above discussion, we have that
for any $\y^+, \y^- \in \O_N$,
\[
|L_G \wt f_{t-s}(\y^+)-L_G \wt f_{t-s}(\y^-)|\le \Eb \sum_{x \in \Zc} 2C \sum_{y \in \Zc \cup I_x}|\y^+(y)-\y^-(y)|.
\]
Therefore applying the above inequality to $\y_s^+, \y_s^-$, we obtain
\begin{align*}
	\E |L_G \wt f_{t-s}(\y_s^+)-L_G \wt f_{t-s}(\y_s^-)|
	&\le 2C \E \Eb \sum_{x \in \Zc} \sum_{y \in \Zc \cup I_x}|\y_s^+(y)-\y_s^-(y)|\\
	&=2C \Eb \sum_{x \in \Zc}\sum_{y \in \Zc \cup I_x}\E|\y_s^+(y)-\y_s^-(y)|,
\end{align*}
where we have used the Fubini theorem in the last equality.
Since $\y_s^+ \ge \y_s^-$ under the monotone coupling,
taking the average $\y_s^+-\y_s^-$ by the translation on $\Z_N$ in the last term, 
namely, plugging
$\y_{\,\cdot\,+z, s}^+-\y_{\,\cdot\,+z, s}^-$ into the above inequality
and taking the arithmetic mean over $z \in \Z_N$,
we obtain
\begin{align*}
2C \Eb \sum_{x \in \Zc}\sum_{y \in \Zc \cup I_x}\E(S(\y_s^+)-S(\y_s^-))
\le 2C|I|(|I|+|I_0|)\E(S(\y_s^+)-S(\y_s^-)),
\end{align*}
and thus
\[
\frac{1}{N}\sum_{x \in \Z_N}\E |L_G \wt f_{t-s}(\y_{\,\cdot\,+x,s}^+)-L_G \wt f_{t-s}(\y_{\,\cdot\,+x,s}^-)|
\le C' \E(S(\y_s^+)-S(\y_s^-)),
\]
for a constant $C' > 0$, 
where $\y_{\,\cdot\,+x,s}^+:=(\y_{s}^+(y+x))_{y \in \Z_N}$ and $\y_{\,\cdot\,+x,s}^-$ is defined similarly.
We recall that $F(\eta)=(1/N)\sum_{x\in\mathbb Z_N}f(\eta_{\cdot+x})$ and note that $(L_GF)(\eta)=\sum_{x\in\mathbb Z_N}(L_Gf)(\eta_{\cdot+x})$.
Hence integrating in $s$ from $0$ to $t$, we have that by \eqref{Eq:integration-by-parts},
\begin{align*}
\Big|\E [F(\y^+_t)-F(\y^-_t)]-\Eb [F(\y^{+z(N^2 t)})-F(\y^{-z(N^2 t)})]\Big|
\le C'\int_0^t \E(S(\y_s^+)-S(\y_s^-))\,ds.
\end{align*}
Note that the constant $C'$ depends only on $\|f\|_\infty$ and $|\supp f|$ as well as $\|c(0, \cdot )\|_\infty$ and $|\supp c(0, \cdot)|$.
This shows the claim.
\qed

The proof of Lemma \ref{Lem:nonlinear} is based on
the following lemma, which we prove in Section \ref{Sec:taming}.

\begin{lemma}\label{Lem:replacement}
For $K$-distinct sites $x_1,\ldots, x_K$ in $\mathbb Z_N$ where $K\ge1$, let $\psi:\Omega_N\mapsto\R$ be a local function of the form
\begin{align*}
\psi(\eta):=\prod_{i=1}^K\eta(x_i), \quad \text{for $\eta\in\Omega_N$}.
\end{align*}
Let $z(t)=\{z_i(t)\}_{i=1, \dots, K}$ for $t\ge 0$ be a $K$-marked SSEP on $\Z_N$ starting from $\{x_i\}_{i=1, \dots, K}$.
If we define
\[
\D_T(\y):=\Eb \psi(\y^{z(T)})-\prod_{i=1}^K \Eb \y(z_i(T)) \quad \text{for $\y \in \O_N$ and $T \ge 0$},
\]
then for all $0<\e<1/28$ there exist constants $C>0$ and $T_0$ such that for all $T_0 \le T \le N^{2+\e}$
the following holds:
for all $\y^+, \y^- \in \O_N$,
\[
|\D_T(\y^+)-\D_T(\y^-)|\le \frac{C}{T^{\e/4}}\sum_{i=1}^K\Eb |\y^+(z_i(T))-\y^-(z_i(T))|+C \exp(-T^{\e/8}),
\]
where the implied constant $C$ depends only on $K$ and $\e$.
\end{lemma}

\proof[Proof of Lemma \ref{Lem:nonlinear}]
Fix an arbitrary $\wt \y=(\y^+, \y^-)$ with $\y^+ \ge \y^-$.
For each $x \in \Z_N$,
let $z(t)=(z_i(t))_{i=1, \dots, K}$ for $t \ge 0$ be a $K$-marked SSEP
starting from $\{x_i+x\}_{i=1, \dots, K}$ 
and $\Eb_{x}$ be the corresponding expectation.
Then for all $T \ge T_0$ and for all configurations $\y^+, \y^- \in \O_N$, we have
\begin{align}\label{Eq:block}
	&\left|\sum_{x \in \Z_N}\Big(\prod_{i=1}^K \Eb_{x} \y^+(z_i(T))-\prod_{i=1}^K \Eb_{x} \y^-(z_i(T))\Big)
	-\sum_{x \in \Z_N}\Big(\F_x(\y^+, T)^K-\F_x(\y^-, T)^K\Big)\right|\nonumber\\
	&\qquad \qquad \qquad \qquad \qquad \qquad 
	\le \frac{C}{\sqrt{T}}\sum_{x \in \Z_N}\sum_{i=1}^K\left|\Eb_{x} \y^+(z_i(T))-\Eb_{x} \y^-(z_i(T))\right|.
\end{align}
Indeed, recall that $\F_{z_i(0)}(\y, T)=\Eb \y(z_i(T))$, and
a ``smoothing effect" of SRW (Lemma \ref{Lem:TV-SRW}) implies that
for $i \neq j$, 
\begin{align}\label{Eq:loc-ave}
&|\F_{z_i(0)}(\y, T)-\F_{z_j(0)}(\y, T)| 
\le 2\|\Pr_{z_i(0)}(z_i(T) \in \,\cdot\,)-\Pr_{z_j(0)}(z_j(T) \in \,\cdot\,)\|_{\TV}\nonumber\\
&\qquad \qquad \le \frac{2C|x_i-x_j|}{\sqrt{T}}\le \frac{2CD}{\sqrt{T}}\quad \text{for all $T \ge T_0$ and $\y \in \O_N$},
\end{align}
where $D:=\max_{i, j=1 ,\dots, K}|x_i-x_j|$.
For a short hand notation, we write
$\F_i^+:=\F_{z_i(0)}(\y^+, T)$ for $i=1, \dots, K$
and similarly for $\F_i^-$.
%where we suppress $x$.
Let us use the following decomposition:
\begin{align}\label{Eq:prod}
\prod_{i=1}^K \F_i^+ - \prod_{i=1}^K \F_i^-
=\(\F_1^+-\F_1^-\)\prod_{i\neq 1}\F_i^+
+\F_1^-\(\F_2^+-\F_2^-\)\prod_{i \neq 1, 2}\F_i^+
+\cdots
+\left(\prod_{i\neq K}\F_i^-\right)\(\F_K^+-\F_K^-\).
\end{align}
Expanding the product in a similar way, we obtain by \eqref{Eq:loc-ave} 
for $1 \le l < K$,
\[
\Big|\prod_{i=i_1, \dots, i_l}\F_i^+ - (\F_j^+)^l\Big| \le \frac{2CDl}{\sqrt{T}}
\le \frac{2CDK}{\sqrt{T}}
\quad \text{for $j, i_1, \dots, i_l =1, \dots, K$},
\]
for all $T \ge T_0$ and the same for $\F_i^-$.
We combine those estimates and see that by triangle inequalities the difference between $\prod_{i=1}^K \F_i^+ - \prod_{i=1}^K \F_i^-$ and
\begin{align*}
\(\F_1^+-\F_1^-\)(\F_1^+)^{K-1}
+\F_2^-\(\F_2^+-\F_2^-\)(\F_2^+)^{K-2}
+\cdots
+(\F_K^-)^{K-1}\(\F_K^+-\F_K^-\)
\end{align*}
is at most
\[
\frac{2CDK}{\sqrt{T}}\sum_{i=1}^K (\F_i^+-\F_i^-).
\]
Note that the second sum of the left-hand side in \eqref{Eq:block} can be written as the sum of the penultimate display by the translation on $\Z_N$.
Therefore summing up over all $x \in \Z_N$, and abbreviating the constant, we obtain \eqref{Eq:block}.

Denote by $\psi_x(\y)=\prod_{i=1, \dots, K}\y(x+x_i)$.
Combining Lemma \ref{Lem:replacement} with \eqref{Eq:block} for any $\e<1/28$, we have that for all $T_0 \le T \le N^{2+\e}$ (this restriction on the range for $T$ is required to apply Lemma \ref{Lem:replacement}),
\begin{align*}
	&\left|\frac{1}{N}\sum_{x \in \Z_N}\Eb_{x} [\psi_x(\y^{+z(T)})- \psi_x(\y^{-z(T)})]
	-\frac{1}{N}\sum_{x \in \Z_N}[\F_x(\y^+, T)^K-\F_x(\y^-, T)^K]\right|\\
	&\le \left|\frac{1}{N}\sum_{x \in \Z_N}\Eb_{x} [\psi_x(\y^{+z(T)})- \psi_x(\y^{-z(T)})]
		-\frac{1}{N}\sum_{x \in \Z_N}\Big(\prod_{i=1}^K \Eb_{x} \y^+(z_i(T))-\prod_{i=1}^K \Eb_{x} \y^-(z_i(T))\Big)\right|\\
		&\qquad \qquad \qquad \qquad \qquad \qquad \qquad +\frac{C}{\sqrt{T}N}\sum_{x \in \Z_N}\sum_{i=1}^K\left|\Eb_{x} \y^+(z_i(T))-\Eb_{x} \y^-(z_i(T))\right|\\
	&\le \frac{C}{T^{\e/4} N}\sum_{x \in \Z_N}\sum_{i=1}^K\Eb_{x} \(\y^+(z_i(T))-\y^-(z_i(T))\)+C e^{-T^{\e/8}}\\
	&=\frac{CK}{T^{\e/4} N}\sum_{x \in \Z_N}\(\F_x(\y^+, T)-\F_x(\y^-, T)\)+C e^{-T^{\e/8}},
\end{align*}
where we have used \eqref{Eq:block} in the first inequality and Lemma \ref{Lem:replacement} and $\e$ is less than $1/2$ in the second inequality.
If we define $F(\y)=(1/N)\sum_{x \in \Z_N}\psi_x(\y)$,
then
by Lemma \ref{Lem:integration-by-parts} for $t=T/N^2$ and for $\wt \y=(\y^+, \y^-)$ where $\y^+ \ge \y^-$,
\begin{align*}
	\Big|\E_{\wt \y} [F(\y^+_{t})-F(\y^-_{t})]-\Eb [F(\y^{+z(T)})-F(\y^{-z(T)})]\Big|
	&\le C\int_0^t \E(S(\y_s^+)-S(\y_s^-))\,ds,
\end{align*}
and plugging $T=tN^2$ in the penultimate display, we conclude that for all $T_0/N^2 \le t \le N^\e$,
\begin{align}\label{Eq:F}
	&\left|\E_{\wt \y}[F(\y_{t}^+)-F(\y_{t}^-)]-\frac{1}{N}\sum_{x \in \Z_N}[\F_x(\y^+, N^2 t)^K-\F_x(\y^-, N^2 t)^K]\right|\nonumber\\
	&\quad \quad \quad \quad \le \frac{CK}{t^{\e/4} N^{1+2\e}}\sum_{x \in \Z_N}\(\F_x(\y^+, N^2 t)-\F_x(\y^-, N^2 t)\)\\
	& \quad \quad \quad \quad \quad \quad \quad \quad +C\int_0^t \E(S(\y_s^+)-S(\y_s^-))\,ds+C \exp(-(N^{2}t)^{\e/8}).\nonumber
\end{align}
Note that for any $s \ge 0$,
\[
\xi_0=S(\y^+)-S(\y^-)=\frac{1}{N}\sum_{x \in \Z_N}\(\F_x(\y^+, s)-\F_x(\y^-, s)\).
\]
Recall that in the general case $F(\y)$ has the form $F(\y)=(1/N)\sum_{x \in \Z_N}\sum_{|I|\neq1}a_I f_I(\y_{\cdot+x})$
and $G(\rho)=\sum_{|I|\neq1}a_I \rho^{|I|}$.
Taking a linear combination and applying to the above inequality \eqref{Eq:F},
we conclude the proof of Lemma \ref{Lem:nonlinear}.
\qed

\section{Replacement lemma: Proof of Lemma \ref{Lem:replacement}}\label{Sec:taming}

For an integer $K \ge 1$ we consider the $K$-marked SSEP
$z(t):=\{z_i(t)\}_{i=1, \dots, K}$ on $\Z_N$,
and the $K$-marked independent SRWs $z^0(t):=\{z^0_i(t)\}_{i=1, \dots, K}$ with rate $1$,
i.e., each particle runs as a SRW having an independent Poisson clock with intensity $1$ (see the explicit generator in Appendix \ref{Sec:A}).
We use the following coupling between these two processes;
the result is known (e.g.\ \cite{DMPbook}) and has been used extensively in \cite{DMFL}.
Since we need to clarify the dependence on every constant involved in the processes and we also need to refer to the construction of the coupling, 
we include the proof of the following claim.

\begin{lemma}\label{Lem:coupling}
	For any integers $N \ge 2$ and $K \ge 1$ with $N \ge K$,
	fix any initial condition $z(0)=z^0(0)=\{x_1, \dots, x_K\}$ in $\Z_N$ where $x_i \neq x_j$ for $i\neq j$, $i, j=1, \dots, K$.
	Then there exists a coupling $\Pb$ between a $K$-marked SSEP $\{z(t)\}_{t \ge 0}$ and $K$-marked independent SRWs $\{z^0(t)\}_{t \ge 0}$ with rate $1$ such that the following holds: for any $\e>0$ there exists $T_0$ such that for all $T_0 \le T \le N^{2+\e}$,
	\[
	\Pb\(\max_{j=1, \dots, K}\max_{0 \le t \le T}|z_j(t)-z_j^0(t)| \ge K T^{\frac{1}{4}+3\e}\) \le 4 K^2\exp\(-T^{\e/4}\),
	\]
	and $z_1(t)=z_1^0(t)$ for all $t \ge 0$ almost surely.
\end{lemma}

In the following, we will use several results on the SRW on $\Z_N$;
all of them are classical --- we include the proofs for the sake of completeness in Appendix \ref{Sec:A} and \ref{Sec:misc}.

\proof[Proof of Lemma \ref{Lem:coupling}]
We construct a coupling from the graphical representation on $\Z_N \times [0, \infty)$.
For each bond $\{x, x+1\}$ in $\Z_N$, we define random exclusion times $\Ec^{x, x+1}=\{\Ec_n^{x, x+1}; n \ge 1\}$ where $\Ec^{x, x+1}$ are independent Poisson point processes with intensity $1$.
Furthermore, for each bond $\{x, x+1\}$, we define a sequence of independent $1/2$-Bernoulli random variables $\{U^{x, x+1}_n; n \ge 1\}$ where $U^{x, x+1}_n \in \{0, 1\}$ for $n \ge 1$ and these sequences are independent of each other and of everything else.
For each $t=\Ec^{x, x+1}_n$, we draw a double arrow between $(x, t)$ and $(x+1, t)$ with mark $U_n^{x, x+1}$.

Given the initial condition $z_i(0)=x_i$ for $i=1, \dots, K$,
the $K$-particles evolve as follows: 
for $t > s \ge 0$, given $z_i(s)$, $i=1, \dots, K$, we put $z_i(t)=z_i(s)$ for all $i=1, \dots, K$ as far as there is no double arrow in time $(s, t]$.
Starting from time $t=0$, for each $i=1, \dots, K$, if there is a double arrow with mark $1$ between $(z_i(t-), t)$ and $(x, t)$, then we let $z_i(t)=x$, and do nothing otherwise.
%Starting from time $t=0$, they exchange the sites $z_i(t-)$ and $z_j(t-)$ so that $z_i(t)=z_j(t-)$ and $z_j(t)=z_i(t-)$ whenever there is a double arrow with mark $1$ between $(z_i(t-), t)$ and $(z_j(t-), t)$, and do nothing otherwise.
Then $\{(z_i(t))_{i=1, \dots, K}\}_{t \ge 0}$ is almost surely defined and has the distribution of $K$-marked SSEP starting from $z_i(0)=x_i$ for $i=1, \dots, K$.

Concerning the $K$-marked independent SRWs, let us define $\{(z^0_i(t))_{i=1, \dots, K}\}_{t \ge 0}$ depending on the number $i=1, \dots, K$.
Suppose that there is a double arrow between $x$ and $x+1$ at time $t$
%$(x, t)$ and $(x+1, t)$ 
for $x \in \Z_N$ and $t \ge 0$.
There are three cases according to the possible configurations on the sites $x, x+1$:
\begin{itemize}
	\item[(1)] There are no particles from $z$ on the sites $x, x+1$.
	In this case we do nothing.
	\item[(2)] There is a single particle from $z$ on the sites $x, x+1$. 
	Suppose that $z_i(t-)=x$ or $x+1$. 
	If the double arrow has the mark $1$,
	then we define 
	\[
	z^0_i(t):=z^0_i(t-)+(z_i(t)-z_i(t-)),
	\]
	and keep everything else, and if the double arrow has the mark $0$, then we do nothing.
	\item[(3)] There are two particles $z_i(t-), z_j(t-)$ on the sites $x, x+1$.
	Assume that $i<j$ and the double arrow has the mark $U \in \{0, 1\}$.
	Then we define
	\[
	U=1 \implies z_i^0(t):=z_i^0(t-)+(z_i(t)-z_i(t-)) \quad \text{and} \quad z_j^0(t):=z_j^0(t-),
	\]
	and
	\[
	U=0 \implies z_j^0(t):=z_j^0(t-)+(z_i(t)-z_j(t-)) \quad \text{and} \quad z_i^0(t):=z_i^0(t-),
	\]
	and we keep everything else.
\end{itemize}
The construction gives $K$-independent SRWs $\{(z^0_i(t))_{i=1, \dots, K}\}_{t \ge 0}$ with rate $1$ (cf.\ Figure 1). %(cf.\ Figure \ref{Fig:coupling})

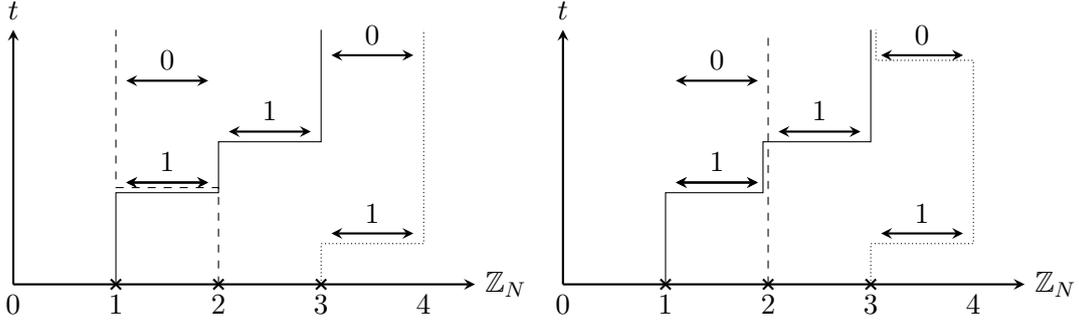
\begin{figure}\label{Fig:coupling}
\centering
\begin{tabular}{cc}
{\begin{tikzpicture}[x=1.5cm,y=1.5cm,scale=0.9]
\coordinate (O1) at (0,2.5) node at (O1) [above] {$t$};
\coordinate (O2) at (4.5,0) node at (O2) [right=0] {$\mathbb{Z}_N$};
\coordinate (A) at (0,0) node at (A) [below] {$0$};
\coordinate (B) at (1,0) node at (B) [below] {$1$};
\coordinate (C) at (2,0) node at (C) [below] {$2$};
\coordinate (D) at (3,0) node at (D) [below] {$3$};
\coordinate (E) at (4,0) node at (E) [below] {$4$};
\coordinate (B1) at (1.5,1.2) node at (B1){$1$};
\coordinate (B2) at (1.5,2.2) node at (B2){$0$};
\coordinate (C1) at (2.5,1.7) node at (C1){$1$};
\coordinate (D1) at (3.5,0.7) node at (D1) {$1$};
\coordinate (D2) at (3.5,2.45) node at (D2) {$0$};
%\fill (B1) circle [radius=2pt];
%\fill (B2) circle [radius=2pt];
%\fill (C1) circle [radius=2pt];
%\fill (D1) circle [radius=2pt];
%\fill (D2) circle [radius=2pt];
\draw[->,>=stealth, thick] (0,0)--(4.5,0);
\draw[->,>=stealth, thick] (0,0)--(0,2.5);
\draw[<->,>=stealth, thick] (1.1,1)--(1.9,1);
\draw[<->,>=stealth, thick] (1.1,1)--(1.9,1);
\draw[<->,>=stealth, thick] (1.1,2)--(1.9,2);
\draw[<->,>=stealth, thick] (2.1,1.5)--(2.9,1.5);
\draw[<->,>=stealth, thick] (3.1,0.5)--(3.9,0.5);
\draw[<->,>=stealth, thick] (3.1,2.25)--(3.9,2.25);
\draw (1,0)--(1,0.9)--(2,0.9)--(2,1.4)--(3,1.4)--(3,2.5);
\draw[dashed] (2,0)--(2,0.95)--(1,0.95)--(1,2.5);
\draw[densely dotted] (3,0)--(3,0.4)--(4,0.4)--(4,2.5);
\draw[thick] (0.95,-0.05)--(1.05,0.05);
\draw[thick] (0.95,0.05)--(1.05,-0.05);
\draw[thick] (1.95,-0.05)--(2.05,0.05);
\draw[thick] (1.95,0.05)--(2.05,-0.05);
\draw[thick] (2.95,-0.05)--(3.05,0.05);
\draw[thick] (2.95,0.05)--(3.05,-0.05);
\end{tikzpicture}}

{\begin{tikzpicture}[x=1.5cm,y=1.5cm,scale=0.9]
\coordinate (O12) at (0,2.5) node at (O12) [above] {$t$};
\coordinate (O22) at (4.5,0) node at (O22) [right=0] {$\mathbb{Z}_N$};
\coordinate (A2) at (0,0) node at (A2) [below] {$0$};
\coordinate (B2) at (1,0) node at (B2) [below] {$1$};
\coordinate (C2) at (2,0) node at (C2) [below] {$2$};
\coordinate (D2) at (3,0) node at (D2) [below] {$3$};
\coordinate (E2) at (4,0) node at (E2) [below] {$4$};
\coordinate (B12) at (1.5,1.2) node at (B12) {$1$};
\coordinate (B22) at (1.5,2.2) node at (B22) {$0$};
\coordinate (C12) at (2.5,1.7) node at (C12) {$1$};
\coordinate (D12) at (3.5,0.7) node at (D12) {$1$};
\coordinate (D22) at (3.5,2.45) node at (D22) {$0$};
%\fill (B12) circle [radius=2pt];
%\fill (B22) circle [radius=2pt];
%\fill (C12) circle [radius=2pt];
%\fill (D12) circle [radius=2pt];
%\fill (D22) circle [radius=2pt];
\draw[->,>=stealth, thick] (0,0)--(4.5,0);
\draw[->,>=stealth, thick] (0,0)--(0,2.5);
\draw[<->,>=stealth, thick] (1.1,1)--(1.9,1);
\draw[<->,>=stealth, thick] (1.1,1)--(1.9,1);
\draw[<->,>=stealth, thick] (1.1,2)--(1.9,2);
\draw[<->,>=stealth, thick] (2.1,1.5)--(2.9,1.5);
\draw[<->,>=stealth, thick] (3.1,0.5)--(3.9,0.5);
\draw[<->,>=stealth, thick] (3.1,2.25)--(3.9,2.25);
\draw (1,0)--(1,0.9)--(1.95,0.9)--(1.95,1.4)--(3,1.4)--(3,2.5);
\draw[dashed] (2,0)--(2,2.5);
\draw[densely dotted] (3,0)--(3,0.4)--(4,0.4)--(4,2.2)--(3.05,2.2)--(3.05,2.5);
\draw[thick] (0.95,-0.05)--(1.05,0.05);
\draw[thick] (0.95,0.05)--(1.05,-0.05);
\draw[thick] (1.95,-0.05)--(2.05,0.05);
\draw[thick] (1.95,0.05)--(2.05,-0.05);
\draw[thick] (2.95,-0.05)--(3.05,0.05);
\draw[thick] (2.95,0.05)--(3.05,-0.05);
\end{tikzpicture}}
\end{tabular}
\caption{Left: $3$-marked SSEP starting from $x_1=1, x_2=2, x_3=3$.
Right: $3$-independent SRWs starting from $x_1=1, x_2=2, x_3=3$.
Three paths (which are solid, dotted, and densely dotted, respectively) correspond to trajectories of three marked particles in the processes, respectively.}
\end{figure}

Let us show that the resulting coupling $\Pb$ is the desired one.
We analyze $z_i(t)-z_i^0(t)$ for $t \ge 0$ and for each $i=1, \dots, K$.
Note that in this construction we have $z_1^0(t)=z_1(t)$ for all $t \ge 0$.
%Namely, the $1$-st independent SRW runs exactly the same as the $1$-st particle in the $K$-marked SEP.
For $j>1$,
$z_j(t)-z_j^0(t)$ can differ only during the period when $|z_j(t)-z_i(t)|=1$ for some $i<j$.
For each $i<j$, denoting by $w_{i, j}(t)$ the difference created in that period up to time $t$,
we have
\[
z_j(t)-z_j^0(t)=\sum_{i<j}w_{i, j}(t) \quad \text{for $j=2, \dots, K$}.
\]
For $i<j$, if we define $\t_\star:=\inf\{t \ge 0 \ : \ |z_i(t)-z_j(t)|=1\}$, then $\t_\star<\infty$ almost surely, and
given $\t_\star$ and $z_i(\t_\star), z_j(\t_\star)$, 
let
\[
T_\star:=\inf\{s \ge 0 \ : \ |z_i(s+\t_\star)-z_j(s+\t_\star)|>1\}.
\]
Then the strong Markov property shows that $T_\star$ conditioned on the $\s$-algebra $\Gc_{\t_\star}$ associated with the filtration $\Gc_t:=\s(\{z(s)\}_{0 \le s \le t})$ for $t \ge 0$ and the stopping time $\t_\star$ has the exponential distribution with rate $1$, and is independent of the random exchange times between $z_i(\t_\star)$ and $z_j(\t_\star)$ since $T_\star$ depends only on the random exchange times after $\t_\star$ on $(x-1, x)$ and $(x+1, x+2)$ for $\{z_i(\t_\star), z_j(\t_\star)\}=\{x, x+1\}$.
Let us define for $t \ge 0$ and $i<j$,
\[
\th_{i, j}(t):=|\{s \in [0, t] \ : \ |z_i(s)-z_j(s)|=1\}|, 
\]
where $|\cdot|$ is the normalized Lebesgue measure.
Conditioned on the random variable $\th_{i, j}(t)$,
each $w_{i, j}(t)$ has the distribution of a SRW with rate $1$ starting from $0$ on $\Z_N$ at time $\th_{i, j}(t)$.
(Note however that for each $j>2$, $\{w_{i, j}(t)\}_{i<j}$ are not independent.)

We observe that for each $i, j$ with $i<j$, the random variable $\th_{i, j}(t)$ is stochastically dominated by
the occupation time of SRW with rate $2$ on the 3-sites $\{-1, 0, 1\}$ on $\Z_N$ starting from $z_0:=z_i(0)-z_j(0)$.
Hence by Lemma \ref{Lem:occupation},
letting $\th(T)$ be the occupation time at $0$ up to time $T$ for the continuous-time SRW on $\Z_N$ with rate $2$ starting from $0$,
we have that
for any $\e>0$, there exists $T_0$ such that for all $T_0 \le T \le N^{2+\e}$ and for each pair $i, j$ with $i<j$,
\begin{equation}\label{Eq:theta}
\Pr\(\th_{i, j}(T) \ge 6 T^{\frac{1}{2}+2\e}\) \le 3\Pr\(\th(T) \ge 2 T^{\frac{1}{2}+2\e}\) \le 3\exp\(-T^{\e/4}\).
\end{equation}
Then, it holds that for all large enough $T$ satisfying that $T \le N^{2+\e}$
for each pair $i<j$,
\begin{align*}
	\Pr\(\max_{0 \le t \le \th_{i,j}(T)} |w_{i,j}(t)| \ge T^{\frac{1}{4}+3\e}\)	
	&\le \Pr\(\max_{0 \le t \le 6 T^{\frac{1}{2}+2\e}} |w_{i,j}(t)| \ge T^{\frac{1}{4}+3\e}\)+\Pr\(\th_{i,j}(T) \ge 6 T^{\frac{1}{2}+2\e}\)\\
	&\le 3\exp\(-\frac{1}{6}T^{2\e}\)+3\exp\(-T^{\e/4}\)
	\le 4\exp\(-T^{\e/4}\),
\end{align*}
where we have used Lemma \ref{Lem:SRW}
with $(6T^{\frac{1}{2}+2\e})^{\frac{1}{2}+\e}<T^{\frac{1}{4}+3\e}$
for $T \ge T_0$ (a re-chosen $T_0$ if necessary)
and \eqref{Eq:theta} in the second inequality.

Noting that for each $j>1$ and for $t, T >0$, 
\begin{align*}
	\Pb\(\max_{0 \le t \le T}|z_j(t)-z_j^0(t)| \ge KT\) \le \sum_{i<j}\Pb\(\max_{0 \le t \le T}|w_{i, j}(t)| \ge T\),
\end{align*}
we have that for any $\e>0$ and for all large enough $T$ with $T \le N^{2+\e}$,
\[
\Pb\(\max_{0 \le t \le T}|z_j(t)-z_j^0(t)| \ge K T^{\frac{1}{4}+3\e}\) \le 4 K\exp\(-T^{\e/4}\).
\] 
Therefore the union bound over $j=1, \dots, K$ implies that for any $\e>0$ there exists $T_0'$ such that for all $T_0' \le T \le N^{2+\e}$,
\[
\Pb\(\max_{j=1, \dots, K}\max_{0 \le t \le T}|z_j(t)-z_j^0(t)| \ge K T^{\frac{1}{4}+3\e}\) \le 4 K^2\exp\(-T^{\e/4}\),
\]
as required.
\qed

In order to prove Lemma \ref{Lem:replacement},
we use a coupling between a $K$-marked SSEP and $K$-marked independent SRWs as it is constructed in Lemma \ref{Lem:coupling}.
A heuristic explanation for the proof of Lemma \ref{Lem:replacement} is as follows.
Let $s_T:=T-T^{\frac{1}{2}+7\e}$. Decomposing the time interval $[0,T]$ as $[0,s_T]$ and $[s_T, T]$,
we construct some events (which will be denoted by $\Gc, \Hc_1$ and $\Hc_2$)
to ensure that $K$-particles in the $K$-marked SSEP are fairly well-separated at time $s_T$ and behave like independent SRWs
during the time range $[s_T,T]$. Since in the first stage $[0, s_T]$ those particles are already fairly separated, they remain sufficiently separated during $[s_T, T]$.
Moreover, the `good' event $\Gc \cap \Hc_1 \cap \Hc_2$ happens with high probability. The construction of these events is a key ingredient
to prove Lemma \ref{Lem:replacement}.

\proof[Proof of Lemma \ref{Lem:replacement}]
Recall that we consider a local function
$\psi(\y)=\prod_{i=1}^K \y(x_i)$ for $K$-distinct sites $x_1, \dots, x_K$ in $\Z_N$.
We consider a coupling $\Pb$ between $\{z(t)\}_{t\ge 0}$ and $K$-marked independent SRWs $\{z^0(t)\}_{t \ge 0}$ with rate $1$ such that $z_i(0)=z_i^0(0)$ for $i=1, \dots, K$ as it is constructed in Lemma \ref{Lem:coupling}.
For all $0<\e<1/28$,
let $s_T:=T-T^{\frac{1}{2}+7\e}$.
First we define the event $\Gc$ where it holds that for all $i=1, \dots, K$,
if $i \neq j$, then
\[
|z_i^0(s_T)-z_j^0(s_T)| \ge T^{\frac{1}{2}-\e}.
\]
Note that $z_1^0(t)=z_1(t)$ holds for all $t \ge 0$ almost surely in the coupling $\Pb$.
By a classical estimate for the SRW (Lemma \ref{Lem:LLT-SRW})
applied to
\[
z_i^0(t)-z_j^0(t) \quad \text{for $i, j=1, \dots, K$ with $i \neq j$ conditioned on $\{z_1^0(t)\}_{0 \le t \le T}$},
\]
there exists $T_0$ such that for all $T_0 \le T \le N^{2+\e}$,
we have almost surely in $\Pb$,
\begin{equation*}
\Pb\(\Gc^c \mid \{\s(z^0_1(t))\}_{0 \le t \le T}\) \le \frac{CK^2}{T^{\e/4}},
\end{equation*}
in particular, since $z^0_1(t)=z_t(t)$ for all $t \ge 0$ almost surely in $\Pb$, we have that
\begin{equation}\label{Eq:Gc}
\Pb\(\Gc^c \mid \s(z_1(T))\) \le \frac{CK^2}{T^{\e/4}} \qquad \text{almost surely in $\Pb$},
\end{equation}
where $\Gc^c$ stands for the complement of the event $\Gc$.
Next we define the event $\Hc_1$ where it holds that for all $i=1, \dots, K$,
\[
|z_i(s_T)-z_i^0(s_T)| \le KT^{\frac{1}{4}+3\e}.
\]
Then Lemma \ref{Lem:coupling} shows that
\begin{equation}\label{Eq:H1}
\Pb\(\Hc_1^c\) \le 4K^2 \exp(-T^{\e/4}) \qquad \text{for all $T_0 \le T \le N^{2+\e}$},
\end{equation}
where we take another $T_0$ if necessary and denote it by the same symbol.
Based on the coupling $\Pb$, we couple $\{z(t)\}_{s_T \le t <T}$ and $K$-marked independent SRWs $\{w^0(t)\}_{t \ge 0}$ with rate $1$
such that $w^0_i(0)=z_i(s_T)$ for each $i=1, \dots, K$ in the same way as in Lemma \ref{Lem:coupling}-----we denote the resulting coupling by the same symbol $\Pb$.
Let us define the event $\Hc_2$ where it holds that for all $i=1, \dots, K$,
\[
|z_i(t+s_T)-w_i^0(t)| \le KT^{\frac{1}{8}+4\e} \quad \text{for all $0 \le t \le T^{\frac{1}{2}+7\e}$}
\]
and for all $i \neq j$,
\[
|w_i^0(t)-w_j^0(t)| \ge \frac{1}{2}T^{\frac{1}{4}+\e} \quad \text{for all $0 \le t \le T^{\frac{1}{2}+7\e}$}.
\]
Lemma \ref{Lem:coupling} and an application of the maximal inequality (the second claim of Lemma \ref{Lem:SRW} with $\lambda=2$ applied to $\{w_i^0(t)-w_j^0(t)\}$ for $0 \le t \le T^{\frac{1}{2}+7\e}$ with $i \neq j$) imply that by the Markov property of $\{z(t)\}_{t \ge 0}$,
\begin{align}\label{Eq:H2}
\Pb\(\Hc_2^c \mid \s(z(s_T), z^0(s_T))\) \1_{\Gc} 
&\le 4K^2 \exp(-T^{\e/8})+3K^2\exp(-\frac{1}{12}T^{\e})\nonumber\\
&\le 5K^2 \exp(-T^{\e/8}) \quad \text{for all $T_0 \le T \le N^{2+\e}$},
\end{align}
where we have used
\[
\(T^{\frac{1}{2}+7\e}\)^{\frac{1}{4}+3\e}=T^{\frac{1}{8}+\frac{13}{4}\e+21 \e^2}<T^{\frac{1}{8}+4\e} \quad \text{for $0<\e<\frac{1}{28}$}, \quad \text{and} \quad T^{\(\frac{1}{2}+7\e\)\(\frac{1}{2}+\e\)} > T^{\frac{1}{4}+\e} \quad \text{for $\e>0$},
\]
and further we choose another $T_0$ depending only on $\e$ if necessary.
Note that given the event $\Gc \cap \Hc_1 \cap \Hc_2$, 
we have by the triangle inequality
\begin{align*}
|z_i(t+s_T)-z_j(t+s_T)|
&\ge |w_i^0(t)-w_j^0(t)|-|z_i(t+s_T)-w_i^0(t)|-|z_j(t+s_T)-w_j^0(t)|\\
&\ge \frac{1}{2}T^{\frac{1}{4}+\e}-2K T^{\frac{1}{8}+4\e}
>1
\end{align*}
for all $i\neq j$, for all $0 \le t \le T-s_T$ and for all $T_0 \le T \le N^{2+\e}$, choosing another $T_0$ depending only on $\e$ and $K$ if necessary,
and thus by the construction of coupling between $\{z(t)\}_{s_T \le t \le T}$ and $\{w^0(t)\}_{t \ge 0}$,
we have that for all $i=1, \dots, K$,
\[
z_i(t+s_T)=w^0_i(t) \quad \text{for all $0 \le t \le T-s_T$}.
\]

Let us decompose
\begin{equation}\label{Eq:decomposition}
\Eb \psi(\y^{z(T)})=\Eb \psi(\y^{z(T)})\1_\Gc+\Eb \psi(\y^{z(T)})\1_{\Gc^c},
\end{equation}
and estimate the second term of \eqref{Eq:decomposition}.
For each $i=1, \dots, K$,
if we define
$g_i(\y):=\prod_{j \neq i}\y(x_j)$ for $\y \in \O_N$, 
then $g_i(\y^{z(T)})=\prod_{j \neq i}\y(z_j(T))$, and
we have
\begin{align*}
\Eb \psi(\y^{z(T)})\1_{\Gc^c}
&=\Eb \left[\Eb[\psi(\y^{z(T)})\1_{\Gc^c}\mid \s(z_i(T))]\right]\\
&=\Eb\left[\y(z_i(T)) \Eb[g_i(\y^{z(T)})\1_{\Gc^c}\mid \s(z_i(T))]\right].
\end{align*}
Letting
$r(i, \y ,T):=\Eb[g_i(\y^{z(T)})\1_{\Gc^c}\mid \s(z_i(T))]$,
by \eqref{Eq:Gc} we obtain  
$|r(i, \y, T)| \le C/T^{\e/4}$
for all $T_0 \le T \le N^{2+\e}$.
Note that for all $\y^+, \y^- \in \O_N$,
\begin{align*}
|\Eb (\y^+(z_i(T))-\y^-(z_i(T)))g_i(\y^{+ z(T)})\1_{\Gc^c}| 
&= |\Eb r(i, \y^+, T) \(\y^+(z_i(T))-\y^-(z_i(T))\)| \nonumber\\
&\le \frac{C}{T^{\e/4}} \Eb |\y^+(z_i(T))-\y^-(z_i(T))|.
\end{align*}
For $\y^+, \y^- \in \O_N$,
using the decomposition
\begin{align*}%\label{Eq:prod}
&\prod_{i=1}^K \y^+(x_i)-\prod_{i=1}^K \y^-(x_i)
=(\y^+(x_1)-\y^-(x_1)) \prod_{j \neq 1} \y^+(x_j)\nonumber\\ 
&+\y^-(x_1)(\y^+(x_2)-\y^-(x_2))\prod_{j \neq 1, 2}\y^+(x_j)+\cdots+\(\prod_{i\neq K}\y^-(x_i)\) (\y^+(x_K)-\y^-(x_K)),
\end{align*}
and permuting the role of $x_i$ for $i=1, \dots, K$,
we obtain
\begin{align}\label{Eq:small-Lip}
\left|\Eb \psi(\y^{+ z(T)})\1_{\Gc^c}-\Eb \psi(\y^{- z(T)})\1_{\Gc^c}\right| \le \frac{C}{T^{\e/4}}\sum_{i=1}^K\Eb |\y^+(z_i(T))-\y^-(z_i(T))|.
\end{align}

Concerning the first term of \eqref{Eq:decomposition}, we have by \eqref{Eq:H1} and \eqref{Eq:H2} for all $T_0 \le T \le N^{2+\e}$,
\[
\Eb \psi(\y^{z(T)})\1_\Gc=\Eb \psi(\y^{z(T)})\1_{\Gc\cap \Hc_1\cap \Hc_2}+O\(e^{-T^{\e/8}}\),
\]
and 
$\Eb \psi(\y^{z(T)})\1_{\Gc \cap \Hc_1\cap \Hc_2}=\Eb\left[\Eb \left[\psi(\y^{z(T)})\1_{\Hc_2} \mid \Fc_{s_T}\right] \1_{\Gc\cap \Hc_1}\right]$.
Given $\Gc\cap \Hc_1$, one has
\begin{align*}
\Eb \left[\psi(\y^{z(T)})\1_{\Hc_2} \mid \Fc_{s_T}\right]
=\Eb \prod_{i=1}^K \y(w^0_i(T-s_T)) \1_{\Hc_2}=\Eb \prod_{i=1}^K \y(w^0_i(T-s_T))+O(e^{-T^{\e/8}}),
\end{align*}
and
\begin{align*}
\Eb \prod_{i=1}^K \y(w^0_i(T-s_T))=\prod_{i=1}^K \Eb \y(w^0_i(T-s_T))=\prod_{i=1}^K \F_{z_i(s_T)} (\y, T-s_T),
\end{align*}
where the first equality follows since $\{w^0(t)\}_{t \ge 0}$ are independent SRWs and we recall $\F_x(\y, t)=\Eb_x \y(w(t))$ a local average of $\y$ by the SRW with rate $1$ around $x$.
Then we use the following estimate on SRW: 
for all $t >0$,
\begin{align*}
\|\Pr_{z_i(s_T)} \(w(t) \in \cdot \)-\Pr_{z^0_i(s_T)}\(w(t) \in \cdot\)\|_{\TV} \le \frac{C|z_i(s_T)-z^0_i(s_T)|}{\sqrt{t}},
\end{align*}
(Lemma \ref{Lem:TV-SRW}),
and we have that given $\Hc_1$
for each $i=1, \dots, K$,
\begin{align*}
|\F_{z_i(s_T)} (\y, T-s_T)-\F_{z^0_i(s_T)}(\y, T-s_T)|
&\le 2\|\Pr_{z_i(s_T)} \(w(T-s_T) \in \cdot \)-\Pr_{z^0_i(s_T)}\(w(T-s_T) \in \cdot\)\|_{\TV}\\
&\le\frac{2C|z_i(s_T)-z^0_i(s_T)|}{\sqrt{T-s_T}}
\le \frac{2CKT^{\frac{1}{4}+3\e}}{\sqrt{T^{\frac{1}{2}+7\e}}}
=\frac{2CK}{T^{\e/2}},
\end{align*}
for all $T_0 \le T \le N^{2+\e}$.
This implies that given $\Hc_1$,
for each $i=1, \dots, K$,
\begin{align}\label{Eq:1}
\prod_{j\neq i} \F_{z_j(s_T)}(\y, T-s_T)=\prod_{j \neq i}\F_{z^0_j(s_T)}(\y, T-s_T)+O\(\frac{1}{T^{\e/2}}\).
\end{align}
For the simplicity of notation, let
\[
\F_i^+:=\F_{z_i(s_T)}(\y^+, T-s_T) \quad \text{and} \quad \F_i^{+0}:=\F_{z_i^0(s_T)}(\y^+, T-s_T),
\]
and similarly $\F_i^-$ and $\F_i^{-0}$.
Since $z_1(t)=z_1^0(t)$ under the coupling, we have $\F_1^\pm=\F_1^{\pm0}$.
Hence by \eqref{Eq:1} we obtain given $\Hc_1$,
\begin{align*}
\(\F_1^+-\F_1^-\)\prod_{j\neq 1}\F_j^+
&=\(\F_1^{+0}-\F_1^{-0}\)\prod_{j\neq 1}\F_j^+\\
&=\(\F_1^{+0}-\F_1^{-0}\)\prod_{j\neq 1}\F_j^{+0}+O\(\frac{1}{T^{\e/2}}\) \(\F_1^{+0}-\F_1^{-0}\).
\end{align*}
Changing the role of $i=1,\ldots,K$, we obtain given $\Hc_1$
\begin{align*}
\left( \prod_{j<i}\Phi_j^- \right) (\Phi_i^+ - \Phi_i^-) \left( \prod_{j>i}\Phi_j^+ \right) =  \left( \prod_{j<i}\Phi_j^{-0}  \right) (\Phi_i^{+0} - \Phi_i^{-0}) \left( \prod_{j>i}\Phi_j^{+0} \right) +O\left( \dfrac{1}{T^{\e/2}} \right) (\Phi_i^{+0} - \Phi_i^{-0}).
\end{align*}
Then using \eqref{Eq:prod} and the same decomposition for $\F_i^{+0}$ and $\F_i^{-0}$ we obtain given $\Hc_1$,
\begin{align*}
\prod_{i=1}^K \F_i^+-\prod_{i=1}^K \F_i^-=\prod_{i=1}^K \F_i^{+0}-\prod_{i=1}^K \F_i^{-0}+O\(\frac{1}{T^{\e/2}}\)\sum_{i=1}^K\(\F_i^{+0}-\F_i^{-0}\).
\end{align*}
Therefore
\begin{align}\label{Eq:2}
\Eb \(\psi(\y^{+z(T)})-\psi(\y^{-z(T)})\)\1_\Gc
&=\Eb\(\prod_{i=1}^K \F_i^{+0}-\prod_{i=1}^K \F_i^{-0}\)\1_{\Gc\cap \Hc_1}\notag\\
&\qquad+O\(\frac{1}{T^{\e/2}}\)\sum_{i=1}^K\Eb\(\F_i^{+0}-\F_i^{-0}\)\1_{\Gc\cap \Hc_1}
+O(e^{-T^{\e/8}}).
\end{align}
By a similar discussion as we did for \eqref{Eq:small-Lip}, we have
\begin{align}\label{Eq:6}
\Eb\left|\prod_{i=1}^K \F_i^{+0}-\prod_{i=1}^K \F_i^{-0}\right|\1_{\Gc^c}
=O\(\frac{1}{T^{\e/4}}\) \sum_{i=1}^K\Eb\left|\F_i^{+0}-\F_i^{-0}\right|.
\end{align}
Finally, we note that
\begin{equation}\label{Eq:3}
\Eb \prod_{i=1}^K\F_{z^0_i(s_T)}(\y, T-s_T)=\prod_{i=1}^K \Eb \F_{z^0_i(s_T)}(\y, T-s_T)=\prod_{i=1}^K \F_{z_i(0)}(\y, T),
\end{equation}
where the first equality holds since $\{z^0(t)\}_{t \ge 0}$ are independent SRWs 
and the second equality follows from the Markov property of each $z^0_i(t)$.
Thus, by \eqref{Eq:decomposition}, \eqref{Eq:small-Lip}, \eqref{Eq:2}, \eqref{Eq:6} and \eqref{Eq:3}, we conclude that
\begin{align*}
|\D_T(\y^+)-\D_T(\y^-)|
&\le \Eb\left|\prod_{i=1}^K \F_i^{+0}-\prod_{i=1}^K \F_i^{-0}\right|\1_{\Gc^c\cup \Hc_1^c}+O\(\frac{1}{T^{\e/2}}\)\sum_{i=1}^K\Eb\left|\F_i^{+0}-\F_i^{-0}\right|\\
&+\frac{C}{T^{\e/4}}\sum_{i=1}^K\Eb |\y^+(z_i(T))-\y^-(z_i(T))|+ O(e^{-T^{\e/8}})\\
&\le \frac{C}{T^{\e/4}}\sum_{i=1}^K\Eb |\y^+(z_i(T))-\y^-(z_i(T))|+O(e^{-T^{\e/8}}),
\end{align*}
for all $T_0 \le T \le N^{2+\e}$, and this concludes the proof of Lemma \ref{Lem:replacement}.
\qed

\appendix

\section{Classical results on SRW}\label{Sec:A}

Fix an integer $N \ge 2$, let $\Z_N:=\Z/N\Z$ be the state space,
where we identify $\Z_N$ with natural representatives $\{0, 1, \dots, N-1\}$ by abusing notations.
We consider the SSEP $\{\y_t\}_{t \ge 0}$ given by the generator
\[
L_E f(\y)=\frac{1}{2}\sum_{x \in \Z_N}(f(\y^{x, x+1})-f(\y)), \qquad \y \in \{0, 1\}^{\Z_N},
\]
where $f: \{0, 1\}^{\Z_N} \to \R$, 
and the continuous-time SRW $\{w(t)\}_{t \ge 0}$ with rate $\lambda >0$ defined by the generator
\[
L_N^{0, \lambda} f(x)=\frac{\lambda}{2} \sum_{y=x\pm 1}\(f(y)-f(x)\), \qquad x \in \Z_N,
\]
where $f:\Z_N \to \R$.
Note that the SRW with rate $\lambda$ runs as follows: a random particle has a Poisson clock with intensity $\lambda$,
and it moves to one of two neighbors with equal probability at each time when the clock rings.

\begin{lemma}\label{Lem:SRW}
	Fix $\lambda>0$.
	Let $\{w(t)\}_{t\ge 0}$ be a continuous-time SRW on $\Z_N$ with rate $\lambda$ and $w(0)=0$.
	Then for any $\e>0$ and for any $T>0$, we have
	\[
	\Pr\(\max_{0 \le t \le T}|w(t)| \ge T^{\frac{1}{2}+\e} \) \le 3\exp\(-\frac{\lambda}{6}T^{2 \e}\).
	\]
	Moreover, if $|w(0)|=\lceil T^{\frac{1}{2}+\e}\rceil<N/2$, then for any $\e>0$ and for any $T >0$, we have
	\[
	\Pr\(\min_{0\le t \le T}|w(t)| \le \frac{1}{2}\lceil T^{\frac{1}{2}+\e}\rceil\) \le 3\exp\(-\frac{\lambda}{24} T^{2 \e}\).
	\]
\end{lemma}

\proof
Let us consider a discrete-time SRW $\{\wt S_t\}_{t=0, 1, \dots}$ on $\Z$ with $\wt S_0=0$.
If we define $S_t:=\wt S_t \mod N$, then $\{S_t\}_{t=0, 1, \dots}$ is a discrete-time SRW on $\Z_N$ with $S_0=0$.
Then the maximal inequality yields for any positive integers $M \ge 1$ and $D \ge 1$,
\[
\Pr\(\max_{0\le k \le M}|S_k| \ge D\) \le \Pr\(\max_{0 \le k \le M}|\wt S_k| \ge D\) \le 2 e^{-\frac{D^2}{2M}}.
\]
For each $T>0$,
the number of jumps $n_T$ in time $[0, T]$ in the continuous-time SRW $\{w(t)\}_{t \ge 0}$ with rate $\lambda$
has the Poisson distribution with intensity $\lambda T$,
and thus we have
\begin{equation}\label{Eq:Lem:occupation}
\Pr\(n_T \ge 3\lambda T\)\le e^{-3\lambda T}\E e^{n_T}=e^{-3\lambda T}e^{(e-1)\lambda T} \le e^{-\lambda T}.
\end{equation}
Therefore we obtain
\begin{align*}
	\Pr\(\max_{0 \le t \le T}|w(t)| \ge D\)	&\le \Pr\( \max_{0 \le k \le n_T}|\wt S_k| \ge D, n_T <3\lambda T\)+\Pr\(n_T \ge 3\lambda T\)\\
	&\le \Pr\( \max_{0 \le k \le \lfloor 3\lambda T\rfloor}|\wt S_k| \ge D\)+e^{-\lambda T}
	\le 2e^{-\frac{D^2}{6\lambda T}}+e^{-\lambda T},
\end{align*}
and letting $D=\lambda T^{\frac{1}{2}+\e}$ for any $T >0$ and for any $\e>0$ yields the first claim.

The second claim follows from the above discussion applying to a continuous-time SRW $w(t)-w(0)$ for $t \ge 0$ starting from $0$;
we omit the details.
\qed

For an integer $K \ge 1$ we consider the $K$-marked SSEP 
$z(t):=\{z_i(t)\}_{i=1, \dots, K}$ on $\Z_N$ running according to the SSEP,
and the $K$-marked independent SRWs $z^0(t):=\{z^0_i(t)\}_{i=1, \dots, K}$ with rate $\lambda$,
where each particle runs as a SRW having an independent Poisson clock with intensity $\lambda$.

We have used the following lemma in
the proof of Lemma \ref{Lem:coupling} in Section \ref{Sec:taming}.

\begin{lemma}\label{Lem:occupation}
	Let $\{w(t)\}_{t \ge 0}$ be a continuous-time SRW on $\Z_N$ with $w(0)=0$ and rate $2$.
	For any $T >0$,
	we define the occupation time at $0$ for $\{w(t)\}_{t\ge 0}$ up to time $T$ by
	\[
	\th(T):=|\{t \in [0, T] \ : \ w(t)=0\}|,
	\]
	where $|A|$ stands for the normalized Lebesgue measure of a measurable set $A \subset \R$.
	Then for any $\e>0$ there exists $T_0$ such that for all $T_0 \le T \le N^{2+\e}$,
	\[
	\Pr\(\th(T) \ge 2 T^{\frac{1}{2}+2\e}\) \le \exp\(-T^{\e/4}\).
	\]
\end{lemma}

\proof
For any positive real $\lambda>0$,
let $\{w(t)\}_{t \ge 0}$ be a continuous-time SRW on $\Z_N$ with $w(0)=0$ and rate $\lambda$
to which we will apply $\lambda=2$.
Since for the continuous-time SRW $\{w(t)\}_{t \ge 0}$ each waiting time for the next jump has the exponential distribution with parameter $\lambda$,
the occupation time $\th(T)$ has the distribution $\sum_{i=1}^{V_T}T_i$ where $V_T$ is the number of visits of $\{w(t)\}_{t \ge 0}$ at $0$ up to time $T$ 
and $T_i$, $i=1, 2, \dots$ are independent exponential random variable with parameter $\lambda$.
We use a discrete-time SRW $\{S_t\}_{t=0, 1, \dots}$ on $\Z_N$ with $S_0=0$ to obtain a bound for $V_T$ and for the occupation time $\th(T)$.

Since the times of jumps for $w(t)$ are given by a Poisson point process on $[0, \infty)$ with intensity $\lambda$,
the number of jumps $n_T$ up to time $T$ satisfies that
%has the Poisson distribution of parameter $\lambda T$.
%Hence 
for all $T \ge 0$, by \eqref{Eq:Lem:occupation},
\begin{equation*}%\label{Eq:Lem:occupation}
\Pr\(n_T \ge 3\lambda T\) \le e^{-\lambda T}.
\end{equation*}

If we denote by $U_n$ the number of visits at $0$ up to time $n$ for discrete-time SRW $\{S_t\}_{t=0, 1, \dots}$,
then by \eqref{Eq:Lem:occupation}, 
\begin{align*}
	\Pr\(V_T \ge M\)	\le \Pr\(V_T \ge M \cap \{n_T < 3 \lambda T\}\)+\Pr\(n_T \ge 3 \lambda T\)
	\le \Pr\(U_{\lfloor 3 \lambda T\rfloor} \ge M\)+e^{-\lambda T}.
\end{align*}
Let $\{\wt S_t\}_{t=0, 1, \dots}$ be a discrete-time SRW on $\Z$ with $\wt S_0=0$
and for $x \in \Z$, let $\wt U_n(x)$ be the number of visits at $x$ up to time $n$ for $\{\wt S_t\}_{t=0, 1, \dots}$.
Noting that $\wt S_t \mod N$ has the same distribution as $S_t$,
we observe that for any positive integer $L>0$,
\begin{equation}\label{Eq:Lem:occupation1}
\Pr\(U_n \ge M\) \le \sum_{x \in (-L, L]\cap N\Z}\Pr\(\wt U_n(x) \ge M \Big\lceil\frac{2L}{N}\Big\rceil^{-1}\)+\Pr\(\max_{1 \le k \le n}|\wt S_k| > L\),
\end{equation}
where we have taken the union bound over $\lceil2L/N\rceil$ points in $(-L, L] \cap N\Z$.
For each $x \in \Z$, it holds that for any integer $m>0$,
\begin{equation}\label{Eq:Lem:occupation2}
\Pr\(\wt U_n(x) \ge m\) \le \Pr\(\wt U_n(0) \ge m\) \le \(1-\frac{1}{2 \sqrt{n}}\)^{m-1} \le 2\exp\(-\frac{m}{2 \sqrt{n}}\),
\end{equation}
for all large enough $n$ independent of $x$; indeed the first inequality holds by the stochastic domination as $\wt S_t$ starts at $0$,
and the second inequality follows since the first return time $R$ to $0$ after $\wt S_t$ leaves $0$ satisfies that
$\Pr\(R > 2n \) \ge (2 \sqrt{n})^{-1}$
for all large enough $n$.
Noting that the maximal inequality yields
\[
\Pr\(\max_{1 \le k \le n}|\wt S_k| \ge L\) \le 2 \exp\(-\frac{L^2}{2n}\),
\]
we have by \eqref{Eq:Lem:occupation1} and \eqref{Eq:Lem:occupation2} 
there exists $n_0$ such that for all integers $M, N$ and $L>0$ and for all large enough $n\ge n_0$,
\[
\Pr\(U_n \ge M\) \le \Big\lceil\frac{2 L}{N}\Big\rceil \exp\(-\frac{M N}{4L \sqrt{n}}\)+2 \exp\(-\frac{L^2}{2n}\).
\]
Summarizing the above estimates, we have that there exists $T_0$  (depending on $n_0$) such that for all $T \ge T_0$
and all integers $M, L>0$,
\[
\Pr\(V_T \ge M\) \le \Big\lceil\frac{2L}{N}\Big\rceil\exp\(-\frac{M N}{4L \sqrt{{\lfloor 3 \lambda T\rfloor}}}\)+2 \exp\(-\frac{L^2}{2{\lfloor 3 \lambda T\rfloor}}\)+e^{-\lambda T}.
\]
Hence fixing the parameter $\lambda=2$, and for any $\e>0$ 
letting $M=\lfloor T^{\frac{1}{2}+2\e}\rfloor$ and $L=\lfloor T^{\frac{1}{2}+\e}\rfloor$,
we obtain
\begin{align*}
\Pr\(V_T \ge T^{\frac{1}{2}+2\e}\) 
\le \frac{2 T^{\frac{1}{2}+\e}}{N} \exp\(-\frac{T^\e N}{4\sqrt{6 T}}\)+2 \exp\(-\frac{T^{2\e}}{12}\)+e^{-2 T}.
\end{align*}
If $T \le N^{2+\e}$, then $T^{1-\e}\le N^2$ and thus 
the first term in the right hand side is at most
\[
2 T^{3\e/2}\exp\(-\frac{T^{\e/2}}{4\sqrt{6}}\) \le \exp\(-T^{\e/3}\),
\]
for all large enough $T$.
Therefore
for any $\e>0$ there exists $T_0$ such that for all $T_0 \le T \le N^{2+\e}$,
\begin{align}\label{Eq:Lem:occupation3}
\Pr\(V_T \ge T^{\frac{1}{2}+2\e}\) 
\le 3\exp\(-T^{\e/3}\).
\end{align}

Returning to the estimate on the occupation time $\th(T)$,
we have that for a sum of independent exponential random variables with parameter $2$,
for any positive integer $V \ge 1$,
\[
\Pr\(\sum_{i=1}^V T_i \ge 2 V\) \le e^{-2V}\prod_{i=1}^V \E e^{T_i}=e^{-2V}\cdot 2^{V} \le e^{-V},
\]
where we have used $\E e^{T_i}=2$ for each $i$.
Therefore combining with \eqref{Eq:Lem:occupation3}, we have that
for any $\e>0$, there exists $T_0'$ such that for all $T_0' \le T \le N^{2+\e}$,
\begin{align*}
	\Pr\(\th(T) \ge 2T^{\frac{1}{2}+2\e}\)	&\le \Pr\(\sum_{i=1}^{V_T}T_i \ge 2 \lfloor T^{\frac{1}{2}+2\e} \rfloor, V_T < T^{\frac{1}{2}+2\e}\)+\Pr\(V_T \ge T^{\frac{1}{2}+2\e}\)\\
	&\le \exp\(-\lfloor T^{\frac{1}{2}+2\e}\rfloor\)+3\exp\(-T^{\e/3}\)
	\le  \exp\(-T^{\e/4}\),
\end{align*}
and this shows the claim.
\qed

\section{Proofs of miscellaneous lemmas}\label{Sec:misc}

\begin{lemma}\label{Lem:TV-SRW}
Fix $\lambda>0$.
There exists a constant $C>0$ such that the following holds:
for $x, y \in \Z_N$, let $\{w_1(t)\}_{t \ge 0}$ and $\{w_2(t)\}_{t \ge 0}$ be continuous-time SRWs on $\Z_N$ with rate $\lambda$ such that $w_1(0)=x$ and $w_2(0)=y$.
Then for all $t >0$,
\[
\|\Pr_x\(w_1(t) \in \cdot \)-\Pr_y\(w_2(t) \in \cdot\)\|_{\TV} \le \frac{C|x-y|}{\sqrt{t}}.
\]
\end{lemma}

\proof
Let us run $\{w_1(t)\}_{t \ge 0}$ and $\{w_2(t)\}_{t \ge 0}$ independently; then
$w_1(t)-w_2(t)$ is a continuous-time SRW on $\Z_N$ with rate $2 \lambda$ starting from $x-y$.
Letting $\t_0$ be the first time $w_1(t)-w_2(t)$ to hit $0$, we have that for any $t \ge 0$,
\begin{equation*}
\|\Pr_x\(w_1(t) \in \cdot \)-\Pr_y\(w_2(t) \in \cdot\)\|_{\TV} \le \Pr\(\t_0 >t\).
\end{equation*}
If we consider a discrete-time SRW $\{S_t\}_{t=0, 1, \dots}$ on $\Z$ starting from $|x-y|$ and let $\t_\ast$ be the first time when it hits $0$,
then there exists a constant $n_0>0$ such that for any $x, y \in \Z$ for $n >n_0$,
\[
\Pr\(\t_\ast>n\) \le \frac{6|x-y|}{\sqrt{n}},
\]
\cite[Theorem 2.17]{LevinPeresBook}.
Since for a continuous-time chain the number of jumps happening in $[0, t]$ has the Poisson distribution with intensity $2 \lambda t$,
and
\begin{align*}
\sum_{n=1}^\infty e^{-2 \lambda}\frac{(2 \lambda t)^n}{n !}\frac{1}{\sqrt{n}}
=\sum_{|n-2\lambda t| \le \lambda t}e^{-2\lambda}\frac{(2 \lambda t)^n}{n !}\frac{1}{\sqrt{n}}+\sum_{|n-2 \lambda t|>\lambda t}e^{-2\lambda}\frac{(2 \lambda t)^n}{n !}\frac{1}{\sqrt{n}}
\le \frac{1}{\sqrt{\lambda t}}+\frac{2}{\lambda t},
\end{align*}
where we have used the Chebyshev inequality in the second term,
we obtain
\[
\Pr\(\t_0 >t\) \le \sum_{n=1}^\infty e^{-2 \lambda}\frac{(2 \lambda t)^n}{n !} \Pr\(\t_\ast >n\)\le \frac{18|x-y|}{\sqrt{\lambda t}},
\]
for $\lambda t \ge n_0$, and upon replacing the constant factor $18/\sqrt{\lambda}$ by a positive constant $C$, we have the claim as stated.
\qed

\begin{lemma}\label{Lem:LLT-SRW}
Fix $\lambda>0$.
Let $\{w(t)\}_{t \ge 0}$ be a continuous-time SRW on $\Z_N$ with rate $\lambda$ such that $w(0)=0$.
Suppose that we have a sequence $s_T$ such that $s_T/T \to 1$ as $T \to \infty$.
Then for any $\e>0$, there exists a constant $C$ such that 
for all large enough $T$ and $N$ with $T \le N^{2+\e}$
and for all $y \in \Z_N$,
\[
\Pr\(|w(s_T) -y| \le T^{\frac{1}{2}-2\e}\) \le \frac{C}{T^{\e/2}}.
\]
\end{lemma}

\proof
Let $\{S_t\}_{t=0, 1, \dots}$ be a discrete-time SRW on $\Z_N$ with $S_0=0$.
For any subset $I$ in $\Z_N$ of size $|I|\le T^{\frac{1}{2}-2\e}$,
let us fix a subset $I_0$ of size $|I|$ in $\{0, 1, \dots, N-1\}$ such that $I=I_0 \mod N$,
and define $\wt I:=[-\lceil t^{1/2+\e} \rceil, \lceil t^{1/2+\e}\rceil ]\cap (I_0+N\Z)$ in $\Z$ for $t \ge 0$.
Taking a SRW $\{\wt S_t\}_{t=0, 1, \dots}$ on $\Z$ such that $S_t=\wt S_t \mod N$,
we apply to $\wt S_t$ on $\Z$ the local limit theorem and the maximal inequality:
\begin{align*}
\Pr(S_{t} \in I) 
& \le \Pr\Big(\wt S_t \in \wt I\Big)+\Pr\Big(\max_{0 \le k \le t}|\wt S_k| \ge t^{1/2+\e}\Big)\\
&\le |I|\frac{t^{\frac{1}{2}+\e}}{N}(1+o(1))\frac{1}{\sqrt{\pi t}}+2 \exp\(-\frac{t^{2\e}}{2}\) \quad \text{as $t \to \infty$},
\end{align*}
where we have used $|\wt I| \le |I| t^{\frac{1}{2}+\e}/N$ in the second inequality.
Note that the number of jumps up to time $s_T$ in a continuous-time SRW is $\lambda s_T$ with an additive error at most $\lambda s_T/2$ with probability at least 
$1-4/\lambda s_T$.
Hence if $s_T/T \to 1$, then for any $\e>0$ and for all large enough $T$ and $N$ with $T \le N^{2+\e}$ (where we recall that $|I| \le T^{\frac{1}{2}-2\e}$),
\begin{align*}
\Pr\(w(s_T) \in I\)&\le T^{\frac{1}{2}-2\e}\frac{T^{\frac{1}{2}+\e}}{N}(1+o(1))\frac{1}{\sqrt{\pi T}}+2 \exp\(-(1-o(1))\frac{T^{2\e}}{2}\)+O\(\frac{1}{T}\)\\
&=O\(\frac{T^{\frac{1}{2}-\e}}{N}\)+O\( \exp\(-\frac{T^{2\e}}{3}\)\)+O\(\frac{1}{T}\)=O\(\frac{1}{T^{\e/2}}\).
\end{align*}
Therefore for any $\e>0$
for all large enough $T$ and $N$ with $T \le N^{2+\e}$ and for any $y \in \Z_N$,
\[
\Pr\(|w(s_T) -y| \le T^{\frac{1}{2}-2\e}\) =O\(\frac{1}{T^{\e/2}}\)
\]
as desired.
\qed

\subsection*{Acknowledgements} 
The authors would like to thank Professors Yasuaki Hiraoka, Shin-ichi Ohta and Tomoyuki Shirai for providing us the opportunity on this collaboration, Hong-Quan Tran and Justin Salez for informing us a flaw in an earlier version of our paper, and anonymous referees for their thorough reading, especially one referee for her or his very careful reviewing, which led to substantial improvement on the presentation.
R.T.\ is supported by JSPS Grant-in-Aid for Scientific Research (C) Grant Number JP20K03602 and JST, ACT-X Grant Number JPMJAX190J, Japan.
K.T.\ is supported by JSPS Grant-in-Aid for Early-Career Scientists Grant Number 18K13426 and 22K13929.

\bibliographystyle{alpha}
\bibliography{GlauberExclusion}

\end{document}